\documentclass[a4paper,12pt]{amsart}
\usepackage{amsmath}
\usepackage{amsthm}
\usepackage{amscd}
\usepackage{amssymb}
\usepackage{amsfonts}
\usepackage{latexsym}
\usepackage[mathscr]{eucal}

\newtheorem{thm}{Theorem}[section]
\newtheorem{prop}[thm]{Proposition}
\newtheorem{lemma}[thm]{Lemma}

\newtheorem{defn}[thm]{Definition}
\newtheorem{exam}{Example}[section]

\newtheorem{rem}{Remark}[section]

\def\supp{\mathop{\rm {supp}}\nolimits}

\def\supp{\hbox{\rm supp}}

\def\A{{\mathcal A}}

\def\O{{\mathcal O}}

\def\Aut{\hbox{\rm Aut}}
\def\K{{\mathcal K}}

\def\F{{\mathcal F}}

\def\A{{\mathcal A}}

\def\C{{\mathcal C}}

\def\bR{{\mathbf R}}

\begin{document}
\title{KMS states on C${}^*$-algebras associated with self-similar sets}
\author{Tsuyoshi Kajiwara}
\address[Tsuyoshi Kajiwara]{Department of Environmental and 
Mathematical Sciences, 
Okayama University, Tsushima, 700-8530,  Japan}      

\author{Yasuo Watatani}
\address[Yasuo Watatani]{Department of Mathematical Sciences, 
Kyushu University, Hakozaki, 
Fukuoka, 812-8581,  Japan}
\maketitle
\begin{abstract}
In this paper, we study KMS states for the gauge actions on
 C${}^*$-algebras associated with self-similar sets whose 
branch points are finite.  
If the self-similar set does not contain any branch point, 
the Hutchinson measure gives the unique KMS state.  But if the 
 self-similar set dose contain a branch point, there sometimes appear 
other KMS states which come from branch points.  For this purpose we 
 construct explicitly a basis for a Hilbert C${}^*$-module associated with 
a self-similar set with finite branch condition.  Using this we get 
condition for a  Borel probability measure on K to be extended to 
a KMS state on the C${}^*$-algebra associated with the original 
self-similar set.  
We classify KMS states for the case of dynamics of unit interval and 
the case of Sierpinski gasket which is related with Complex dynamical 
system.  KMS states for these examples are unique and 
given by the Hutchinson measure if $\beta$ is equal to $\log N$, where 
$N$ is the number of contractions.  They are expressed as convex             
 combinations of KMS states given by measures supported on the orbit 
of the branched points if $\beta > \log N$.
\end{abstract}

\section{Introduction} 
There exist many interactions between reversible topological dynamical 
systems and their C${}^*$-algebras through the crossed product 
construction by groups of homeomorphisms on compact Hausdorff spaces.  
On the other hand we have many interesting examples of irreversible 
dynamical systems of continuous maps like a tent map on the unit 
interval and rational functions on the Riemann sphere.  
They are often branched covering maps or expansive maps  
on compact metric spaces and their inverse branches sometimes consist of 
proper contractions.  Thus a family of proper contractions on a compact 
metric space which is self-similar with respect to these contractions 
give a irreversible dynamical system in some sense.  
Although C${}^*$-algebras of groupoids by Renault \cite{R1} and \cite{R2}
are also useful for irreversible systems, 
we study C${}^*$-algebras of bimodules by Pimsner \cite{Pi} to include 
singular points.  In our study we show that there exists a relation 
between the orbit structure of branched points for irreversible 
dynamical systems and the structure of KMS states for the gauge 
actions on their $C^*$-algebras. 
In this paper we study irreversible systems defined by family of
proper contractions on a self-similar set.  
In \cite{KW2} we introduced a C${}^*$-algebra 
$\O_{\gamma}= \O_{\gamma}(K)$ associated with a system 
$\gamma = (\gamma_1,\dots,\gamma_N)$ of contractions on 
a self-similar set $K$.  
We explicitly determine the KMS states in the case of contractions
corresponding to a tent map on the unit interval and the rational function
$R(z) = \frac{z^3-\frac{16}{27}}{z}$ on the Julia set $J_R$, which is 
homeomorphic to the  Sierpinski gasket. Since Sierpinski gasket 
contains three branched points and their inverse orbits  
by contractions fall in fixed points, we have that for any 
$\beta > \log 3$, the set of $\beta$-KMS sates is homeomorphic 
to a two-dimensional simplex spanned by the three vertices 
corresponding to the three branched points. 
\par 
We recall that Olsen-Pedersen \cite{OP} showed that KMS state with
inverse temperature $\beta$ (i.e.,$\beta$-KMS state) on Cuntz algebra
$\O_n$ with respect to the gauge action exists if and only if 
$\beta = \log n$ and that $\log n$-KMS state is unique.  
Evans \cite{Ev} extended their result for quasi-free actions.  
Enomoto-Fujii-Watatani \cite{EFW} studied the gauge 
action on Cuntz-Krieger algebras $\O_A$, and show that the KMS state is
unique and its inverse temperature is the logarithm of
the Perron-Frobenius eigenvalue of $A$ for a irreducible matrix $A$.  
Exel-Laca \cite{EL} studied KMS states on partial crossed product 
C${}^*$-algebras by free groups and classify KMS 
states on Cuntz-Krieger algebras associated with infinite matrix. Exel 
studied KMS states more in \cite{Ex1} and \cite{Ex2}.  
More generally KMS states on Cuntz-Pimsner algebras are studied by 
Pinzari-Watatani-Yonetani \cite{PWY}, Kerr and Pinzari \cite{KP}
and Laca-Neshveyev \cite{LN}.   
Kumjian and Renault \cite{KR} investigated KMS states on groupoid
C${}^*$-algebras associated with expansive maps.  
In many cases the logarithm of inverse temperature of a KMS 
state is equal to the entropy of the corresponding dynamical 
systems like sofic shifts \cite{MWY}, \cite{PWY}.  
But except the value of entropy, we have been unable to catch 
any information of the structure of the dynamical systems from 
the property of KMS states.  
The aim of the paper is to get a information on the structure of 
branched points from the structure of KMS states.  
\par
We introduced a C${}^*$-algebra $\O_R$ associated with a rational
function $R$ in \cite{KW1} and show that if the Julia set does not 
contain any branched point, then the Lyubich measure gives the unique 
$\log \deg R$-KMS state for the gauge action. 
\par
In this paper, we study  KMS states on C${}^*$-algebras associated with
self-similar sets whose branched points are finite.  If the self-similar 
set does not contain any branched point, then the Hutchinson measure 
gives the unique KMS state for the gauge action.  
But if the self-similar set does contain a branched point, 
there sometimes appear another KMS state.  
To study it, we need to  construct a concrete countable basis for
Hilbert C${}^*$-bimodules.  
First we recall the definition and fundamental results for basis for
Hilbert C${}^*$-module, which will be used without saying explicitly.
The fact that each basis automatically converges unconditionally
is important and used in several occasions.    
Next, we characterize a KMS state on $\O_X$ in terms of its restriction 
to the coefficient algebra $A$, where $A$ is a unital C${}^*$-algebra
and $X$ is a countably generated Hilbert C${}^*$-bimodule over $A$.  
This comes from a general theorem in Laca-Neshveyev \cite{LN}.  But they 
extend traces on $A$ to Toeplitz algebra in some specific class first, 
extend general traces by perturbation of action and take weak limits.  
We here provide a simple and direct proof using the properties of 
countable basis for clear understanding of extension from traces on $A$
to the fixed point algebra of $\O_X$ by the gauge action.  
Next, we explicitly construct a basis called a patched basis,
for Hilbert C${}^*$-modules constructed from self-similar sets
which satisfy the finite branch condition.
Using a patched basis, we express the condition that tracial 
states on $A=\text{C}(K)$ extend to KMS states on $\O_{\gamma}$ without
using basis.  We obtain some Ruell-Perron-Frobenius like operator
concerning the condition that measure on $K$ is extended to KMS states.
Last, we classify KMS states for some specific examples.  
We treat dynamics on unit interval.  It is shown that there exists a 
unique $\log N$-KMS states.  This state is of infinite type defined 
in \cite{LN}.  When there exists a branched value in $K$, another 
type of KMS states appear.  For each branched point $y$,  
there exists a KMS states which is expressed as a countable sum of 
Dirac measures.  These are of finite type in \cite{LN}.  The KMS 
states of these C${}^*$-algebra are expressed by convex combinations 
of them.  When a family of contraction is the section of an map $h$, 
the minimum value of the logarithm of inverse temperature of KMS states
is shown to be the entropy of $h$.  We do similar classification of 
KMS states for the C${}^*$-algebra associated with Sierpinski gasket 
introduced in \cite{KW2}.  
\par
The content of this paper is as follows:  
In section 2, we recall several definitions and fundamental facts.  
In section 3, we give a characterization of KMS states on 
Cuntz-Pimsner algebra using a countable basis.  
In section 4, we construct basis with the finite branch condition and 
provide a characterization of KMS states in terms of measures on 
the self-similar set. 
In section 5, we present classification results of KMS states for specific
examples.  
\par
The method presented in this paper to classify KMS states on the  C${}^*$
algebra expressed as a Cuntz-Pimsner algebra with an abelian 
C${}^*$-algebra as a coefficient algebra is applicable to many other 
cases, and we hope that our method shed light on the role of 
branched points in Cuntz-Pimsner algebra constructed from 
correspondences with branches.  
We shall study KMS states on the C${}^*$-algebra associated 
with complex dynamical systems in the forthcoming papers.   
\par
The authors are partially supported by Grants-in-Aid 
for Scientific Research 15540207 and 14340050 from Japan Society for
the Promotion of Science.

\section{Self-similar sets and Hilbert C${}^*$ bimodules} 
Let $(K,d)$ be a compact metric space.  
\begin{defn}
 A continuous map $\gamma$ on $K$ is called a proper contraction if
 there exists constants $0  < c_1 \le c_2 < 1$ such that 
\[
 c_1 d(y_1,y_2) \le d(\gamma(y_1),\gamma(y_2)) \le c_2 d(y_1,y_2)
 \qquad y_1,\,\,y_2 \in K.  
\]
\end{defn}
Let $N$ be an integer greater than $1$, and let 
$\gamma = (\gamma_1,\dots,\gamma_N)$ be a family of proper 
contractions on $(K,d)$.  We use notations such as  
$\gamma(x) = \bigcup_{j=1}^N \{ \gamma_j(x)\}$ and $\gamma^{-1}(x)
 = \bigcup_{j=1}^N \{ \gamma_j^{-1}(x)\}$.  
\begin{defn}
 $K$ is called self-similar with respect to $\gamma$ when
 $K = \bigcup_{i=1}^N \gamma_i(K)$.  
\end{defn}

\begin{lemma}
If $K$ is self similar with respect to $\gamma$, $K$ has no
isolated point.  
\label{lem:isol}
\end{lemma}
\begin{proof}
We fix $x \in K$.  Then for each $n$ there exists $(j_1,\dots,j_n)$
such that $x \in \gamma_{j_1}\cdots \gamma_{j_n}(K)$ because $K$ is
self-similar with respect to $\gamma$.  
Since $\gamma_j$'s are proper contraction, the diameter of 
$\gamma_{j_1}\cdots \gamma_{j_n}(K)$, which is homeomorphic to $K$,
 tends to zero as $n \to \infty$.  This shows that $x \in K$ is not 
isolated in $K$.  
\end{proof}

We need additional technical conditions.  
\begin{defn}
 We say that $\gamma = ( \gamma_1,\dots,\gamma_N )$ satisfies the open 
set condition if there exists a non empty open subset $V$ of $K$ 
such that $ \bigcup_{i=1}^N 
 \gamma_i(V) \subset V $ and $\gamma_j(V) \cap \gamma_{j'}(V) 
= \phi$ for $j \ne j'$.  
\end{defn}

Define subsets $B(\gamma)$, $C(\gamma)$ and $\tilde{C}(\gamma)$ of 
$K$ by 
\begin{align*} 
 B(\gamma) & = \{ x \in K | x =\gamma_j(y)=\gamma_{j'}(y) 
\text{ for some } y 
 \in K \text{ and } j \ne j' \} \\
 C(\gamma) & =  \{ y \in K | \gamma_j(y)=\gamma_{j'}(y) \text{ for } 
 j \ne j' \} \\
 \tilde{C}(\gamma) & =  \bigcup_{j=1}^{N} \gamma_j^{-1}(B(\gamma)).  
\end{align*}
We call a point in $B(\gamma)$ a branched point,  and a point in
$C(\gamma)$ a branched value.  
\begin{defn}
We say that $\gamma$ satisfies the finite branch condition if
 $C(\gamma)$ is a finite set.  
\end{defn}

We define a branched index $e(x,y)$ when 
$x \in \gamma(y)$ by 
\[
 e(x,y) = \verb!#!\{ j \in \{1,\dots,N \} | \gamma_j(y) = x \}
\]
For $x \in K$ we define $I(x)$ by 
\[
 I(x) = \{ j \in \{1,\dots,N \} | \text{ there exists } y \in K
 \text{ such that } x = \gamma_j(y) \}.  
\]
We also use the following notation: For $y \in K$, 
\[
 O(y)  = \bigcup_{n=0}^{\infty}
 \{\gamma_{j_1}\cdots \gamma_{j_n}(y) | (j_1,\cdots,j_n) \in
  \{ 1,\cdots,N \}^n \}, 
\]
with a convention $\gamma_{j_1}\cdots \gamma_{j_n}(y)=y$ for $n=0$.  
We call $O(y)$ the orbit of $y$.  

\begin{exam}
 Let $K=[0,1]$, $\gamma_1(y) = \frac{1}{2}y$, and $\gamma_2(y) = 1 - \frac{1}{2}y$.  
Then $K$ is self-similar with respect to $\gamma = (\gamma_1,\gamma_2)$.
 $\gamma$ satisfies the open set condition.  In this example,
 we can take $V=(0,1)$.  We have $B(\gamma) =\{ \frac{1}{2}\}$ and 
 $C(\gamma) = \{1\}$.  We refer this example as the case of tent map.  
\label{exam:tent}
\end{exam}

\begin{exam}
 Let $K = [0,1]$, $\gamma_1(y)=\frac{1}{2}y$ and 
$\gamma_2(y) = \frac{1}{2}(y+1)$.  Then $K$ is self-similar with respect 
to $\gamma$.  In this case, $B(\gamma)=\phi$ and $C(\gamma)=\phi$.  
\label{exam:Cuntz}
\end{exam}

We recall Cuntz-Pimsner algebras \cite{Pi}.  
Let $A$ be a C${}^*$-algebra and $X$ be a Hilbert right $A$-module.  
We denote by $L(X)$ the algebra of the adjointable bounded operators 
on $X$.  For $\xi$, $\eta \in X$, the "rank one" operator $\theta _{\xi,\eta}$
is defined by $\theta _{\xi,\eta}(\zeta) = \xi(\eta|\zeta)$
for $\zeta \in X$. The closure of the linear span of rank one 
operators is denoted by $K(X)$.   We say that 
$X$ is a Hilbert bimodule over $A$ if $X$ is a Hilbert right  $A$-
module with a *-homomorphism $\phi : A \rightarrow L(X)$.  We always assume 
that $X$ is full and $\phi$ is injective. 
   Let $F(X) = \bigoplus _{n=0}^{\infty} X^{\otimes n}$
be the full Fock module of $X$ with a convention $X^{\otimes 0} = A$. 
 For $\xi \in X$, the creation operator $T_{\xi} \in L(F(X))$ is defined by 
\[
T_{\xi}(a) =  \xi a  \qquad \text{and } \ 
T_{\xi}(\xi _1 \otimes \dots \otimes \xi _n) = \xi \otimes 
\xi _1 \otimes \dots \otimes \xi _n .
\]
We define $i_{F(X)}: A \rightarrow L(F(X))$ by 
$$
i_{F(X)}(a)(b) = ab \qquad \text{and } \ 
i_{F(X)}(a)(\xi _1 \otimes \dots \otimes \xi _n) = \phi (a)
\xi _1 \otimes \dots \otimes \xi _n 
$$
for $a,b \in A$.  The Cuntz-Toeplitz algebra ${\mathcal T}_X$ 
is the C${}^*$-algebra on $F(X)$ generated by $i_{F(X)}(a)$
with $a \in A$ and $T_{\xi}$ with $\xi \in X$.  
Let $j_K | K(X) \rightarrow {\mathcal T}_X$ be the homomorphism 
defined by $j_K(\theta _{\xi,\eta}) = T_{\xi}T_{\eta}^*$. 
We consider the ideal $I_X := \phi ^{-1}(K(X))$ of $A$. 
Let ${\mathcal J}_X$ be the ideal of ${\mathcal T}_X$ generated 
by $\{ i_{F(X)}(a) - (j_K \circ \phi)(a) ; a \in I_X\}$.  Then 
the Cuntz-Pimsner algebra ${\mathcal O}_X$ is defined as 
the quotient ${\mathcal T}_X/{\mathcal J}_X$ . 
Let $\pi : {\mathcal T}_X \rightarrow {\mathcal O}_X$ be the 
quotient map.  Put $S_{\xi} = \pi (T_{\xi})$ and 
$i(a) = \pi (i_{F(X)}(a))$. Let
$i_K : K(X) \rightarrow {\mathcal O}_X$ be the homomorphism 
defined by $i_K(\theta _{\xi,\eta}) = S_{\xi}S_{\eta}^*$. Then 
$\pi((j_K \circ \phi)(a)) = (i_K \circ \phi)(a)$ for $a \in I_X$.   
We note that  the Cuntz-Pimsner algebra ${\mathcal O}_X$ is 
the universal C${}^*$-algebra generated by $i(a)$ with $a \in A$ and 
$S_{\xi}$ with $\xi \in X$  satisfying that 
$i(a)S_{\xi} = S_{\phi (a)\xi}$, $S_{\xi}i(a) = S_{\xi a}$, 
$S_{\xi}^*S_{\eta} = i((\xi | \eta)_A)$ for $a \in A$, 
$\xi, \eta \in X$ and $i(a) = (i_K \circ \phi)(a)$ for $a \in I_X$.
We usually identify $i(a)$ with $a$ in $A$.  
We also identify $S_{\xi}$ with $\xi \in X$ and simply write $\xi$
instead of $S_{\xi}$.  
We denote by ${\mathcal O}_X^{alg}$ the $\ ^*$-algebra generated
algebraically by $A$  and $X$. There exists an action 
$\alpha : {\mathbb R} \rightarrow \Aut \ {\mathcal O}_X$
with $\alpha_t(\xi) = e^{it}\xi$, which is called the gauge action.
Since we assume that $\phi: A \rightarrow L(X)$ is 
isometric, there is an embedding $\phi _n : L(X^{\otimes n})
 \rightarrow L(X^{\otimes n+1})$ with $\phi _n(T) = 
T \otimes id_X$ for $T \in L(X^{\otimes n})$ with a convention 
$\phi _0 = \phi : A \rightarrow L(X)$.  
There exists an isometric map $j_K^{(n)}$ from $K(X^{\otimes n})$ to 
$\O_{X}$ such that $j_K^{(n)}(\theta_{\xi_1 \otimes \cdots \otimes
\xi_n,\eta_1 \otimes \cdots 
\otimes \eta_n}) = \xi_1 \cdots \xi_n {\eta_n}^* \cdots
{\eta_1}^*$.  We also identify $K(X^{\otimes n})$ and its image in $\O_X$.  
\par
We put 
\[
 \tilde{\F}_n = \phi_{n-1}\circ \cdots \circ \phi_1 \circ \phi (A)
   + \phi_{n-1} \circ \cdots \circ \phi_1(K(X)) + \cdots +
   K(X^{\otimes n}).  
\]
with a convention $\tilde{\F}_0=A$.  
We denote by $\F_n$ the C${}^*$-algebra generated by $A$, $K(X)$,
$\cdots$, and $K(X^{\otimes n})$ in $\O_X$ with a convention $\F_0 = A$.  
Then there exists a family of isomorphisms
 $\{ \Psi_n\}_{n=0}^{\infty}$ between $\tilde{\F}_n$ and $\F_n$, which
 is compatible with two filtrations $\tilde{\F}_{i} \subset \tilde{\F}_{i+1}$
 and $\F_i \subset \F_{i+1}$ \cite{Pi} and \cite{FMR}.  The algebra
 $\F_X = \overline{ \bigcup_{n=0}^{\infty} \F_n}$ coincides with the
 fixed point algebra under the gauge action $\alpha$.  We denote by $E$
 the conditional expectation from $\O_X$ to $\F_X$ given by $\alpha$.  
\par
Let $(K,d)$ be a compact metric space, $\gamma
=(\gamma_1,\cdots,\gamma_N)$ be a system of proper contraction.  
We put $\C_{\gamma}=\bigcup_{j=1}^N \{ (\gamma_j(y),y) | y \in K\}$, 
$A={\rm C}(K)$ and $X_{\gamma}={\rm C}(\C_{\gamma})$.  
For $a$, $b \in A$ and $f$, $g \in X$, we
define two $A$ actions and right $A$-inner product on $X$ as follows:
\[
  (a\cdot f \cdot b)(x,y)  =  a(x)f(x,y)b(y) \qquad 
 (f|g)_A(y)  =  \sum_{j=1}^N
 \overline{f(\gamma_j(y),y)}g(\gamma_i(y),y).  
\] 
Then $X_{\gamma}$ is a full Hilbert $A$-module.  If we put $\phi(a)f=a\cdot f$,
then $\phi$ is an isometric *-homomorphism from $A$ to $L(X)$.  Then
$X_{\gamma}$ is a full Hilbert C${}^*$-bimodule over $A$.  

\begin{defn}
 The Cuntz-Pimsner algebra $\O_{X_{\gamma}}$ constructed from
 $X_{\gamma}$ is denoted by $\O_{\gamma}=\O_{\gamma}(K)$.  
\end{defn}

\begin{thm} \cite{KW2} When $\gamma$ satisfies the open set condition,
 $\O_{\gamma}$  is simple and purely infinite.  
\end{thm}

\begin{lemma} If $\gamma$ satisfies the finite branch 
 condition, we have $I_X = \{a \in A\, | \, a(y)=0 \quad \text{for } 
y \in B(\gamma) \}$.  
\end{lemma} 
\begin{proof}
 This lemma is proved in Proposition 2.4 \cite{KW2} under the open set
 condition.  But the open set condition is only used to take a sequence 
 $(x_n,y_n) \in {\C_{\gamma}}$ such that  
$x_n \notin B(\gamma)$, $x_n \to c$ and 
 $y_n \to d$ for $(c,d)$ with $c \in B(\gamma)$ and $d \in C(\gamma)$.  
 When $K$ is self similar with respect to $\gamma$ and 
 $\gamma$ satisfies the finite branch condition, by Lemma \ref{lem:isol} 
these sequence always
 exists.  
\end{proof}
This lemma also follows from Theorem 3.11 in \cite{MT}.  
\par
For $x \in X$, we write $\|x\|_2 = \| (x|x)_A\|^{1/2}$.  We note that
if $ |f_n(\gamma_j(y),y) - f(\gamma_j(y),y)|  \to 0$ uniformly with
respect to $y$ for every $j$, then we have $\|f_n - f\|_2 \to 0$.  
\par
We recall bases for Hilbert C${}^*$-modules \cite{KPW}.  
In the following, we assume that $A$ is $\sigma$-unital and $X$ is
countably generated.  
A family $(u_{\lambda})_{\lambda \in \Lambda}$ in $X$ indexed by a
set $\Lambda$, which is countable or finite, is called a basis for $X$ if for
every $\varepsilon > 0$ there exists a finite subset $F_0 \subset
\Lambda$ such that for every finite subset $F$ such that $F \supset F_0$
we have 
\[
 \left\| \sum_{k \in F} u_k(u_k|f)_A - f \right\| < \varepsilon, 
\]
that is, a net $\{\sum_{k \in F} u_k(u_k|f)_A | F \subset \Lambda 
\text{ finite subset }\}$ converges to $f$ with respect to $\|
\cdot \|$.  We write this as 
\[
  \sum_{k \in \Lambda} u_k(u_k|f)_A = f.  
\]
We note that for countably infinite $\Lambda$ 
an index set $(u_k)_{k \in \Lambda}$ is a basis if and only 
if for any numbering of $\Lambda$ we have 
\[
 \sum_{k=1}^{\infty}u_k(u_k|f)_A = f
\]
in norm.  This means that the series converges to $f$ unconditionally.  
\par 
We assume that a countable sequence $\{u_k\}_{k=1}^{\infty}$ in $X$ satisfies
\[
 \sum_{k=1}^{\infty} u_k(u_k|f)_A = f
\]
for $f \in X$.  Then by Proposition 1.2 in \cite{KPW}, 
An indexed family $(u_k)_{k \in {\bf N}}$ is automatically a basis for $X$, 
where ${\bf N}$ denotes the set of natural numbers. 
\par
We have the following lemmas.  

\begin{lemma}
 Let $(u^i_k)_{k \in \Lambda^i}$, $(i=1,\dots,n)$ be families in
 $X$ indexed by countable sets $\Lambda^i$.  
Put $\Lambda = \bigcup_{i=1}^n \{(i,k) | k \in \Lambda^i\}$.  
We assume that for every $f \in X$, 
\[
 f =  \sum_{i=1}^n \sum_{k \in \Lambda^i} u^i_k (u^i_k|f)_A.   
\]
Then the indexed family $(u^i_{k} )_{(i,k) \in \Lambda} $ is a basis
 for $X$.   \label{lem:gather}
\end{lemma}
\begin{lemma}
 Let $(u_k)_{k \in \Lambda}$ be a basis for $X$.  
Then $(u_{k_1}\otimes \cdots
  \otimes u_{k_n})_{(k_1,\dots,k_n) \in \Lambda^n}$ is a basis for
 $X^{\otimes n}$.  \label{lem:tensor}
\end{lemma}

\section{KMS states on Pimsner algebras}

Let ${\A}$ be a C${}^*$-algebra with a one parameter automorphism group 
 $\alpha$.  A state $\varphi$ on $\A$ is called a $\beta$-KMS state with
respect to $\alpha $ if 
\[
 \varphi(a \alpha_{i\beta}(b)) = \varphi(b a)
\]
holds for $a \in \A$ and $b \in \A_{a}$, where $\A_{a}$ denotes the set
of entire analytic elements for $\alpha$ in $\A$.  We refer for the
definition and fundamental matters of KMS states to Bratteli-Robinson
\cite{BR}. 
\par
Let $I$ be an ideal of a C${}^*$-algebra $B$.  Let $\psi$ be a positive
linear functional of on $I$.  The natural extension $\tilde{\psi}$ of
$\psi$ to $B$ is given by $\psi(b) = \lim_{\lambda} \psi(b e_{\lambda})$ 
for $b\in B$, where $(e_{\lambda})_{\lambda}$ is an approximate unit in $I$.  
We need the following general lemma.  

\begin{lemma}(Proposition 12.5 \cite{EL}) Let $B$ be a unital 
 C${}^*$-algebra.  Suppose $B=A+I$ where $A$ is a C${}^*$-subalgerba
 containing $1$ and $I$ is a closed two sided ideal.  Let $\phi$ be a
 state on $A$ and $\psi$ be a positive linear functional on $I$.
 We denote by $\tilde{\psi}$ the natural extension of $\psi$ to $B$.  
Then if 
(1) $\phi \ge \tilde{\psi}$ on $A^+$ and (2) $\phi = \psi$ on $A \cap I$, 
there exists a state $\rho$ on $B$ such that $\rho|_A = \phi$ and 
 $\rho|_I=\psi$.  Moreover, such a state $\rho$ is unique.  
\label{lem:ext}
\end{lemma}

Let $A$ be a unital C${}^*$-algebra, $X$ be a countably
generated full Hilbert $A$-module and 
$\phi$ be an injective *-homomorphism from $A$ to $L(X)$.  Let $\O_X$ be
a Cuntz-Pimsner algebra constructed from $ X $.  We also need the following
lemma.  

\begin{lemma}(Lemma 4.2(2) \cite{FMR}) We have
 $\F^{(n-1)} \cap K(X^{\otimes n})
 = K(X^{\otimes n-1}) \cap K(X^{\otimes n})$ in $\O_X$.  
\label{lem:inter}
\end{lemma}

The following proposition follows from the general theorem in \cite{LN}.
 But since the extension procedure in \cite{LN} of traces 
on $A$ to $\O_X$ is not straightforward using perturbation, we give a simple 
and direct proof within  $\O_X$, which is a natural extension of 
Lemma 3.2 in \cite{PWY}.  
We fix a basis $(u_k)_{k \in \Lambda}$.  We assume that $\Lambda$ is
 countably infinite, and admit $u_k=0$ for some $k \in \Lambda$.
 Let $\{u_k\}_{k=1}^{\infty}$ denote the 
sequence obtained from $(u_k)_{k \in \Lambda}$ by some numbering of
$\Lambda$.  

\begin{prop}
The restriction of a $\beta$-KMS state $\varphi$ on $\O_X$ to $A$ is a tracial
 state $\tau$ on $A$ satisfying the following conditions: 
\begin{align*}
  \sum_{k=1}^{\infty} \tau((u_k|au_k)_A) & =  \lambda \tau (a)  
\quad (\forall a \in I_X) \quad (1)  \\
 \sum_{k =1}^{\infty} \tau((u_k|au_k)_A) & \le  \lambda \tau (a) 
\quad (\forall a \in A^{+}) \quad (2) \label{lem:KMS1}
\end{align*} 
for $\lambda = e^{\beta}$.  
\par
For a tracial state $\tau$ on $A$ satisfying (1) and (2), we can
 construct a $\beta$-KMS state $\varphi$ on $\O_X$ whose restriction to $A$
 coincides with $\tau$.  Moreover, such an extension is unique.  
\end{prop}
\begin{proof}

Let $\varphi$ be a $\beta$-KMS state on $\O_{X}$.  By the condition
 of $\beta$-KMS state, we have 
\[
 \sum_{k=1}^n \varphi(u_k^* a u_k) = \lambda \sum_{k=1}^n 
\varphi(a u_i u_i^*)
\]
where $a \in A$.  Since $\sum_{k=1}^n u_k u_k^* \le I$ 
and $\varphi|_{\O_{X}^{(0)}}$ is a trace, we have 

\[
 0 \le  \lambda \sum_{k=1}^n \varphi(a u_k u_k^*) 
   = \lambda \varphi \left( a^{1/2} \left( \sum_{k=1}^n u_k u_k^* \right)
  a^{1/2} \right)  \le \lambda  \varphi(a)
\]
for $a \in A^{+}$.  Then we have, 
\[
 \sum_{k=1}^{\infty} \varphi(u_k^* a u_k) \le  \lambda  \varphi(a).  
\]
Let $a \in I_X \subset \K(X)$.  Since
 $(\sum_{k=1}^{n}\theta_{u_k,u_k})_{n=1}^{\infty}$ is an approximate
 unit in $K(X)$, we have 
\[
  \sum_{k=1}^{\infty} \varphi(u_k^* a u_k)
 =  \lambda  \lim_{n \to \infty}\varphi(a \sum_{k=1}^n \theta_{u_k,u_k})  
 =  \lambda  \varphi(a).  
\] 
Then we conclude that (1) and (2) hold.  
\par
We take a tracial state $\tau$  on $A$ satisfying 
 the condition (1) and (2).  
Let $F_n$ be a finite subset of ${\bf N}^n$.  We define finite sets
 $F_{p-1}$, $\cdots$, $F_{1}$ by $F_{p-1} = \{(k_1,\dots,k_{p-1}) 
\in {\bf N}^{p-1} | (k_1, \dots , k_{p-1},k_p) \in F_p\}$ inductively.
 $(2 \le p \le n)$.   Then using the condition (2) repeatedly,  we have 
\begin{align*}
 &   \lambda^{-n}\sum_{(k_1,\dots,k_n) \in F_n} \tau( (u_{k_1}\otimes
  \cdots \otimes u_{k_n}| 
 u_{k_1}\otimes \cdots \otimes u_{k_n} )_A)  \\
  = & \lambda^{-n}
  \sum_{(k_1,\dots,k_n)\in F_n} \tau (u_{k_n}^* \cdots u_{k_1}^*u_{k_1}
  \cdots u_{u_n}) \\
  \le & \lambda^{-n} \sum_{k_n=1}^{\infty} \sum_{(k_1,\dots,k_{n-1}) \in F_{n-1}}
   \tau ( u_{k_n}^*u_{k_{n-1}}^* \cdots u_{k_1}^*u_{k_1} 
   \cdots u_{k_{n-1}} u_{k_n} )  \\
  \le & \lambda^{-n+1}  \sum_{(k_1,\dots,k_{n-1}) \in F_{n-1}}
  \tau (u_{k_{n-1}}^* \cdots u_{k_1}^*u_{k_1}\cdots u_{u_{n-1}}) \\
  \le & \lambda^{-p} \sum_{(k_1,\dots,k_p) \in F_{p}}
  \tau (u_{k_{p-1}}^* \cdots u_{k_1}^*u_{k_1} \cdots u_{k_{p-1}}) \\
  \le & 1.  
\end{align*}
Since $x \le \|x\| I$ for $x \in L(X^{\otimes n})^+$, we have 
\[
 \lambda^{-n} \sum_{(k_1,\dots,k_n) \in F_n}
 \tau( (u_{k_1}\otimes \cdots \otimes u_{k_n}|
 xu_{k_1}\otimes \cdots \otimes u_{k_n} )_A)  \le \|x\|.  
\]
This shows that $\sum_{(k_1,\dots,k_n) \in \Lambda^n}
 \tau( (u_{k_1}\otimes \cdots \otimes u_{k_n}|
 xu_{k_1}\otimes \cdots \otimes u_{k_n} )_A) $ converges
 unconditionally for each $x \in L(X^{\otimes n})$.  
Then we can define a bounded positive linear functional
 $\sigma^n$ on $L(X^{\otimes n})$ by 
\[ \sigma^{n}(x) =  \lambda^{-n} \sum_{(k_1,\dots,k_n) \in {\bf N}^n}
 \tau( (u_{k_1} \otimes \cdots \otimes u_{k_n}|
 xu_{k_1}\otimes \cdots \otimes u_{k_n} )_A).  
\]
We put $\tau^n = \sigma^n |_{K(X^{\otimes n})}$.  
For $x_1 \otimes \cdots \otimes x_n$, $y_1 \otimes \cdots \otimes y_n
 \in X^{\otimes n}$, we have 
\begin{align*}
 &  \tau^n ( \theta_{x_1 \otimes \cdots \otimes x_n,y_1 \otimes \cdots
 \otimes y_n} )    \\
= & \lambda^{-n}  \sum_{(k_1,\dots,k_n) \in {\bf N}^n}
   \tau ( (u_{k_1}\otimes \cdots \otimes u_{k_n} | \theta_{x_1 
  \otimes \cdots \otimes x_n,y_1 \otimes \cdots \otimes y_n} (u_{k_1} 
  \otimes \cdots \otimes u_{k_n})  )_A) \\
  = & \lambda^{-n} \sum_{(k_1,\cdots,k_n \in) {\bf N}^n} \tau  
((u_{k_1}\otimes \cdots \otimes u_{k_n} | x_1 \otimes
 \cdots \otimes x_n (y_1 \otimes \cdots \otimes y_n|u_{k_1}\otimes
 \cdots \otimes u_{k_n})_A )_A) \\
  = & \lambda^{-n} \sum_{(k_1,\cdots,k_n) \in {\bf N}^n } 
  \tau ((u_{k_1}\otimes \cdots \otimes u_{k_n}|x_1 
  \otimes \cdots \otimes x_n)_A (y_1 \otimes \cdots \otimes y_n|
  u_{k_1}\otimes \cdots \otimes u_{k_n})_A  )\\
  = & \lambda^{-n} \sum_{(k_1,\cdots,k_n ) \in {\bf N}^n } \tau 
 ( (y_1 \otimes \cdots \otimes y_n| u_{k_1}\otimes \cdots \otimes
 u_{k_n})_A 
(u_{k_1}\otimes \cdots \otimes u_{k_n}|x_1 \otimes \cdots \otimes x_n)_A   )\\
  = & \lambda^{-n}\sum_{(k_1,\cdots,k_n) \in {\bf N}^n } 
  \tau ( (u_{k_1}\otimes \cdots \otimes u_{k_n}
  (u_{k_1}\otimes \cdots \otimes u_{k_n}|y_1 \otimes \cdots 
  \otimes y_n)_A | x_1 \otimes \cdots \otimes x_n)_A) \\
  = & \lambda^{-n} \tau ( \sum_{(k_1,\cdots,k_n) \in {\bf N}^n}
 (u_{k_1}\otimes \cdots 
  \otimes u_{k_n} (u_{k_1}\otimes \cdots \otimes u_{k_n}|y_1 \otimes
 \cdots \otimes y_n)_A | x_1 \otimes \cdots \otimes x_n)_A) \\
  = & \lambda^{-n}\tau((y_1 \otimes \cdots \otimes y_n |x_1 \otimes
 \cdots \otimes x_n)_A), 
\end{align*}
because an indexed set $( u_{k_1} \otimes \cdots
 \otimes u_{k_n} )_{(k_1,\dots,k_n) \in \Lambda^n}$ is a basis for
 $X^{\otimes n}$ by Lemma \ref{lem:tensor}.     
\par
Let $F$ be a finite subset of $\Lambda^n$.  We put
 $e_F = \sum_{\{ k_1,\cdots,k_n \}\in F}
 \theta_{u_{k_1}\otimes \cdots \otimes u_{k_n}, u_{k_1}\otimes \cdots
 \otimes u_{k_n}} $.  Then $\{e_F \}_{F \text{ finite subset of }\Lambda^n }$ 
 is an approximate unit in $K(X^{\otimes n})$.  
The we have 
\begin{align*}
 \tilde{\tau}^n(x) & = \lim_{F  } \tau^n (xe_F ) \\
                   & = \lim_{F } \tau^n (\sum_{(k_1,\dots,k_n) \in F}
               \theta_{xu_{k_1} \otimes \cdots \otimes u_{k_n}, u_{k_1}
               \otimes \cdots \otimes u_{k_n} })  \\
               & = \lim_{F }\tau( (u_{k_1}\otimes \cdots \otimes u_{k_n} 
               |xu_{k_1}\otimes \cdots \otimes u_{k_n})_A) \\
                 & = \sigma^n(x).  
\end{align*}
Thus the natural extension $\tilde{\tau}^n$ to $L(X^{\otimes n})$ is given
 by $\sigma^n$. 
Since $\tilde{\F}^{(n-1)} \subset L(X^{\otimes n}) = M(K(X^{\otimes n}))$, the
 natural extension $\tilde{\tau}^n$ of $\tau^n$ to $\F^{(n-1)}$ is given by
 $\sigma^n \circ \Psi_n^{-1}$.  
\par
 For each $n \ge 0$, we define states $\omega^n$ on  $\F^{(n)}$ 
 which extend $\tau$ such that $\omega^{n+1}|_{\F^{(n)}} = \omega^n$
 and $\omega^n|_{K(X^{\otimes n})} = \tau^n$.  
 First, we put $\omega^0= \tau$.  Then $\omega^0$ is a 
 state on $\F^{(0)} = A$.  We assume $n \ge 1$, and assume that there 
 exist states $\omega^i$ on $\F^{(i)}$ for $0 \le i \le n$ such that 
 $\omega^i|_{K(X^{\otimes i})} = \tau^i$ for $1 \le i \le n$ and 
 $\tilde{\tau}^{i} \le \omega^{i}$ on $\F^{(i-1)}$ for $1 \le i \le n$.  
 We have $\F^{(n+1)} = \F^{(n)} + K(X^{\otimes n+1})$, and  by Lemma 
 \ref{lem:inter} we have $\F^{(n)} \cap K(X^{\otimes n+1}) 
= K(X^{\otimes n}) 
 \cap K(X^{\otimes n+1})$.  Let $x \in K(X^{\otimes n}) 
\cap K(X^{\otimes n+1})$.  Using 
 $u_{i_n}^* \cdots u_{i_1} x u_{i_1} \cdots u_{i_n} \in A \cap K(X) = I_X$
 and the condition (1), we have 
\begin{align*}
 \tau^{n+1}(x) & = \lambda^{-n-1} \sum_{(i_1,\cdots,i_n,i_{n+1}) \in
   {\bf N}^{n+1}}   \tau(u_{i_{n+1}}^* u_{i_{i_n}}^*\cdots u_{i_1}^* 
      x u_{i_1} \cdots u_{i_n} u_{i_{n+1}}  ) \\
               & = \lambda^{-n}  \sum_{(i_1,\cdots,i_n ) \in {\bf N}^n}
                \tau( u_{i_{i_n}}^*\cdots u_{i_1}^* x u_{i_1} \cdots
               u_{i_n}  ) \\
               & = \tau^{n}(x) = \omega^n(x).  
 \end{align*}
For $x \in \F^{(n)}$ we write $x = y + z$ where $y \in \F^{(n-1)}$ and 
 $z \in K(X^{\otimes n})$.  By the assumption of induction, we have 
 $\tilde{\tau}^{n}(y^*y) \le \omega^n(y^*y)$.  Since 
 $y^*z + z^*y + z^*z \in K(X^{\otimes n})$, we have 
 $\tau^n(y^*z + z^*y + z^*z) = \omega^n(y^*z + z^*y + z^*z)$.  
 We note that for $x \in {\F^{(n)}}^+$ we have 
\begin{align*}
 \tilde{\tau}^{n+1}(x) & = \lambda^{-n-1}
          \sum_{(i_1,\dots,i_n,i_{n+1}) \in { {\bf N}^{n+1} }}
         \tau( u_{i_{n+1}}^* u_{i_{n}}^* \cdots u_{i_1}^* x
         u_{i_1}\cdots 
            u_{i_{n}} u_{ i_{n+1} } )\\ 
              & \le \lambda^{-n} \sum_{ (i_1,\dots,i_n) 
   \in {{\bf N}^{n+1}}}
         \tau( u_{i_{n}}^* \cdots u_{i_1}^* x u_{i_1}\cdots
            u_{i_{n}} )  \\
           & =  \tilde{\tau}^{n}(x).  
\end{align*}
We used  the fact that $u_{i_{n}}^* \cdots u_{i_1}^* x u_{i_1}\cdots 
u_{i_{n}} \in A^+$  because $x \in {\F^{(n)}}^+$ and the condition (2).  
Then, we have 
\begin{align*}
\tilde{\tau}^{n+1}(x^*x) & = \tilde{\tau}^{n+1}( (y+z)^*(y+z) ) \\
            & \le \tilde{\tau}^{n} ( (y+z)^*(y+z) ) \\
            & = \tilde{\tau}^{n}(y^*y + y^*z + z^*y + z^*z) \\
            & = \tilde{\tau}^{n}(y^*y) + \tau^{n}(y^*z + z^*y + z^*z) \\
            & = \tilde{\tau}^{n}(y^*y) + \omega^{n}(y^*z + z^*y + z^*z) \\
            & \le \omega^{n}(y^*y) + \omega^{n}(y^*z + z^*y + z^*z) \\
            & = \omega^{n}(y^*y + y^*z + z^*y + z^*z ) \\
            & = \omega^{n}(x^*x).    
\end{align*}
By Lemma {\ref{lem:ext}}, there exists a state $\omega^{n+1}$ on
 $\F^{(n+1)} = \F^{(n)} + K(X^{\otimes n+1})$ such that
 $\omega^{n+1}|_{\F^{(n)}} = \omega^n$, and $\omega^{n+1}|_{K(X^{\otimes n+1})}
 = \tau^{n+1}$.  Moreover we have $\tilde{\tau}^{n+1}|_{\F^{(n)}}
 \le \omega^{(n+1)}$.  
\par
By induction, we have states $\omega^n$ on $\F^{(n)}$ for all $n \ge 0$ 
 such that $\omega^{n+1}|_{\F^{(n)}} = \omega^n$, and
 $\omega^n|_{K(X^{\otimes n})} = \tau^{n}$.  Putting
 $\omega(x) = \omega^{n}(x)$ for $x \in \F^{(n)}$, we define $\omega$ on 
$\bigcup_{i=0}^{\infty}\F^{(i)}$.  Since $\omega^n$ is a state on
 $\F^{(n)}$, we can extend $\omega$ to a state on the closure $\F_X$.  
\par
We note that $\omega$ does not depend on the choice of a 
 basis $( u_k )_{k \in \Lambda}$ on $X$ because the values of $\omega$ on
 $K(X^{\otimes n})$ are expressed without basis.  
\par
Let $m$ be an integer.  We show $\omega (x y^*) = 
\lambda^m \omega (y^* x)$ for $x$, $y \in \O_{X}^{(m)}$.  
It is sufficient to prove this equality for nonnegative $m$,  and when 
 $x$ and $y$ are in the form 
\[
x = z_1 \dots z_m x_1 \cdots x_p y_p^* \cdots y_1^*,  \quad
y = z_1' \cdots z_m' x_1' \cdots {x_q'} {y_q'}^* \cdots {y_q'}^*
\]
where $z_i$, $z_i'$, $x_j$, $x_j'$, $y_k$, $y_k' \in X$.  Put 
\[
a = (y_p^* \cdots y_1^*) ({y_1'} \cdots y_p') \in A, \quad 
\tilde{a} = ({x_p'}^* \cdots {x_1'}^* {z_m'}^* \cdots {z_1'}^*)( z_1 \cdots
 z_m x_1 \cdots x_p) \in A.  
\]
We may assume that $p \le q$.  Another case is similar.  
We have 
\begin{align*}
   & \omega(xy^*)   \\
   = & \omega(z_1 \dots z_m x_1 \cdots x_p y_p^* \cdots y_1^* 
  {y_1'} \cdots {y_q'} {x_q'}^* \cdots {x_1'}^* {z_m'}^* \cdots {z_1'}^* )\\
   = & \omega(z_1 \dots z_m x_1 \cdots x_p (y_p^* \cdots y_1^*
  y_1' \cdots y_p') 
   y_{p+1}' \cdots y_q'{x_q'}^* \cdots {x_1'}^* {z_m'}^* \cdots {z_1'}^*) \\
  = & \omega ( z_1 \dots z_m x_1 \cdots x_p a {y_{p+1}'} \cdots {y_q'}
     {x_q'}^* \cdots {x_1'}^* {z_m'}^* \cdots {z_1'}^*) \\
   = & \lambda^{-(q+m)} \tau ({x_q'}^* \cdots {x_1'}^* {z_m'}^* \cdots {z_1'}^*  
     z_1 \dots z_m x_1 \cdots x_p a {y_{p+1}'} \cdots {y_q'}) \\
   = &  \lambda^{-(q+m)} 
   \tau ( {x_q'}^* \cdots {x_1'}^* {z_m'}^* \cdots {z_1'}^*  
     z_1 \dots z_m x_1 \cdots x_p(y_p^* \cdots y_1^* {y_1'} \cdots y_p')
     {y_{p+1}'} \cdots {y_q'}) \\
  = &  \lambda^{-m}\lambda^{-q}
  \tau ( {x_q'}^* \cdots {x_{p+1}'}^* ({x_p'}^* \cdots {x_1'}^* {z_m'}^*
  \cdots {z_1'}^* z_1 \cdots z_m x_1 \cdots  x_p) y_p^* \cdots
  y_1^* {y_1'} \cdots {y_q'}  ) \\
  = &  \lambda^{-m}\lambda^{-q}
  \tau ( {x_q'}^* \cdots {x_{p+1}'}^* \tilde{a}^* y_p^* \cdots y_1^*
  {y_1'} \cdots {y_q'}) \\
  = & \lambda^{-m}
  \omega ( {y_1'} \cdots {y_q'}
  {x_q'}^* \cdots {x_{p+1}'}^* \tilde{a}^* y_p^* \cdots y_1^*) \\
  = &  \lambda^{-m} 
 \omega ( {y_1'} \cdots {y_q'}  {x_q'}^* \cdots {x_{p+1}'}^* 
({x_p'}^* \cdots {x_1'}^* {z_m'}^* \cdots {z_1'}^* z_1 \cdots z_m
 x_1 \cdots  x_p) y_p^*  \cdots y_1^* ) \\
 = & \lambda^{-m}
 \omega (  {y_1'} \cdots {y_q'} {x_q'}^* \cdots {x_1'}^* {z_m'}^* \cdots {z_1'}^*
 z_1 \cdots z_m x_1 \cdots x_p y_p^* \cdots y_1^*  ) \\
  = & \lambda^{-m} \omega(y^*x).  
\end{align*}
When $m=0$, this shows that $\omega$ is a trace.  When $m=1$,
 this shows that $\omega$ can be extended to a $\beta$-KMS state
 $\varphi$  on $\O_{X}$ by $\varphi = \omega\circ E$.  
\end{proof}

If $\Lambda$ is finite, we have $I_X=A$ and the condition (2) is
unnecessary.  In this case, this proposition is Lemma 3.2 of \cite{PWY}.   
If (2) holds, $\sum_{k=1}^{\infty} \tau((u_k|au_k)_A)$ dose not depend
on the numbering of $\Lambda$ and we can write 
$\sum_{k \in \Lambda} \tau(( u_k |a u_k)_A)$.

\section{KMS states on a C${}^*$-algebra associated with a self-similar set}

Let $(K,d)$ be a compact metric space and $\gamma = (\gamma_1,\cdots,
\gamma_N)$ be a set of proper contractions on $K$.  
In the following, we assume that $\gamma$ satisfies the finite branch
condition.  
We construct a basis for Hilbert C${}^*$-module $X_{\gamma}$
over $A$.  
\par 
For this purpose we consider the following situation.  Let $K_1$ and
$K_2$ be compact metric spaces, $n$ be an integer and $\gamma_i$
($i=1,\cdots,n$) be proper contractions from $K_1$ to $K_2$.  
For $n \ge 2$,  we assume that there exists a $c \in K_1$ such that
$\gamma_1(c)= \cdots =\gamma_n(c)$ and for $\gamma_i(y)$'s 
are different $y \ne c$.  We put 
$\C=\{ (\gamma_i(y),y) | y \in K_1,\,\, i=1,\dots,n \}$, 
$A={\rm C}(K_1)$ and $X = {\rm C}(\C)$.  Then $X$ is a right Hilbert 
$A$ module.  We say that such a module is of $n$-branch class.  
We construct a basis for Hilbert C${}^*$-module of $n$-branch class.  
If $n=1$, we put $\Lambda = \{1\}$ and $u_1(x,y)=1$.  
Then $(u_k)_{k \in \Lambda}$ is a basis for $X$.  
\par
Assume that $n \ge 2$.  We fix a positive $P$.
We define a family of a functions $r_i(x)$ on $[0,\infty )$ for $i \ge 1$ by 
 \[
 r_i(x)  = 
 \begin{cases}
  & 1 \qquad \frac{P}{i} \le x \\
  & \left(\frac{i}{2P}\right) x-1 \qquad \frac{P}{2i}
  \le x \le \frac{P}{i} \\
  & 0 \qquad 0 \le x \le \frac{P}{2i}
 \end{cases}
\]
and with a convention $r_0(x)=0$.  For $i \ge 0$, $r_i(x)$ is a non decreasing
function and $r_i(x) \le r_{i+1}(x)$ for every $x$.  We put
$v_i(x) = (r_i(x)-r_{i-1}(x))^{1/2}$ for $i \ge 1$.  
Let $\delta > 0$, then there exists an $i_{\delta} > 0$ such that
$\sum_{\tilde{i}=1}^{i_{\delta}}v_{\tilde{i}}(\delta)^2 = 1$.  
For $x \ge \delta$ and $i \ge i_{\delta}$ we have $v_i(x)=0$, and then 
we have $\sum_{\tilde{i}=1}^{i}v_{\tilde{i}}(x)^2 = 1$.  
\par 
Let $\omega =e^{2 \pi \sqrt{\,-1}/n}$.  Since $\gamma_j(y)$
($j=1,\cdots,n$) are different for $y \ne c$ and $v_i(0)=0$, 
we can do the following definition.  
For $k \ge 1$, we define a family of continuous functions $u_k$ in 
$X$ as follows:  
\begin{align*}
 u_1(x,y) & = \frac{1}{\sqrt{\,n}} \\
 u_{1+(n-1)(i-1) + l}(\gamma_j(y),y)  & = \frac{1}{\sqrt{\,n}}
 \,\omega^{lj}\,  v_i(d(y,c)), 
\end{align*}
where $i \ge 1$ and $1 \le l \le n-1$.  
\par
For $y$ with $d(y,c) \ge \delta$ and $k$
with $k \ge 1 + (n-1)(i-1)$ where $i \ge i_{\delta}$, we have
$u_k(\gamma_j(y),y)=0$ for each $j$.  Let $M \ge 1+(n-1)(i_{\delta}-1)$, 
and put $f_M = \sum_{k=1}^M u_k(u_k|f)_A$.   For each $j$, we have 
\begin{align*}
 f_M( \gamma_j(y),y) & = \sum_{k=1}^M u_k(\gamma_{j}(y),y)
  \sum_{\tilde{j}=1}^n \overline{u_k(\gamma_{\tilde{j}} (y),y)} 
  f(\gamma_{\tilde{j}}(y),y) \\ 
 & = \frac{1}{n}  \sum_{{\tilde{j}}=1}^n  f(\gamma_{\tilde{j}}(y),y) + 
   \frac{1}{n}\sum_{\tilde{i}=1}^{i_{\delta}} \sum_{{\tilde{j}}=1}^n \sum_{l=1}^{n-1} 
    \omega^{l(j - {\tilde{j}})} v_{\tilde{i}}(d(y,c))^2 
   f(\gamma_{\tilde{j}}(y),y) \\
 & = \frac{1}{n}  \sum_{{\tilde{j}}=1}^n   f(\gamma_{\tilde{j}}(y),y)
   + \frac{(n-1)}{n}\sum_{\tilde{i}=1}^{i_{\delta}}
 v_{\tilde{i}}(d(y,c))^2 f(\gamma_{j}(y),y) \\
 &    - \frac{1}{n} \sum_{ \tilde{j}=1,({\tilde{j}} \ne j)}^n   
  \sum_{\tilde{i}=1}^{i_{\delta}} v_{\tilde{i}}(d(y,c))^2
  f(\gamma_{\tilde{j}}(y),y) \\
   & = f(\gamma_{j}(y),y).  
\end{align*} 
We used $\sum_{l=1}^{n-1}\omega^{ l (\tilde{j}-j) }$ is $-1$ for
$\tilde{j} \ne j$ and $n-1$ for $\tilde{j} = j$, 
and $\sum_{\tilde{i}=1}^{i_{\delta}} v_{\tilde{i}}(d(y,c))^2=1$ for $y$ 
with $\delta \le d(y,c)$.  
\par
We take arbitrary small $\varepsilon'>0$.  There exists a $\delta > 0$ 
satisfying the following: If $d(y,c) < \delta$, then we have 
\[
 |f(\gamma_j(y),y) - f(b,c)| < \varepsilon'
\]
for each $j$.  We write as $M = 1 + (n-1)i + l$ with $i \ge 1$ and
$1 \le l \le n-1$.  For $1 \le j \le n$ we have 
\begin{align*}
 & | f_M(\gamma_j(y),y)-f(\gamma_j(y),y)| \\
 = & | f_M ( \gamma_j(y),y)-f(b,c) | + | f( \gamma_j(y),y) - f(b,c) | \\
 \le & | f_M ( \gamma_j (y),y) - f(b,c) | + \varepsilon'.  
\end{align*}

We estimate  $|f_M(\gamma_{j}(y),y)-f(b,c)|$.  
\begin{align*}
 |f_M(\gamma_{j}(y),y)-f(b,c)| 
 \le  & \left| \frac{1}{n} \sum_{{\tilde{j}}=1}^{n} (f(\gamma_{\tilde{j}}(y),y)) -
 f(b,c) \right| \\ 
 &  + \left| \frac{1}{n} \sum_{p=1}^{n-1}\sum_{{\tilde{j}}=1}^n
  \omega^{p({\tilde{j}}-j)} \sum_{\tilde{i}=1}^{i} v_{\tilde{i}}(d(y,c))^2 
 (f(\gamma_{\tilde{j}}(y),y) - f(b,c)) \right|   \\
 &  + \left| \frac{1}{n} \sum_{p=1}^{n-1} \sum_{{\tilde{j}}=1}^n 
 \omega^{p({\tilde{j}}-j)} \sum_{\tilde{i}=1}^{i} v_{\tilde{i}}(d(y,c))^2 f(b,c) \right| \\
 & + \left| \frac{1}{n} \sum_{p=1}^{l} \sum_{{\tilde{j}}=1}^n
  \omega^{p({\tilde{j}}-j)} v_{i+1}(d(y,c))^2 
 (f(\gamma_{\tilde{j}}(y),y)-f(b,c))\right|  \\
 & + \left| \frac{1}{n} \sum_{p=1}^{l} \sum_{{\tilde{j}}=1}^n
   \omega^{p({\tilde{j}}-j)} v_{i+1}(d(y,b))^2 f(b,c)\}\right| \\
 \le & \frac{1}{n} n \varepsilon' + \frac{1}{n} n (n-1)
 \left( \sum_{\tilde{i}=1}^{i} v_{\tilde{i}}(d(y,b))^2\right) \varepsilon'
  + \frac{1}{n} l n v_{i+1}(d(y,b)) \varepsilon' \\
 \le & \varepsilon' + (n-1) \left( \sum_{\tilde{i}=1}^{i+1} v_{\tilde{i}}(d(y,c))^2\right)
 \varepsilon' \\
 \le & n \varepsilon'.  
\end{align*}

We used $\sum_{\tilde{j}=1}^{n}\omega^{p(\tilde{j}-j)} = 0$ for $1 \le p 
\le n-1$ and $\sum_{\tilde{i}=1}^{i+1}v_{\tilde{i}}(d(y,c))^2 \le 1$.  
We take arbitrary small $\varepsilon$, and take $\varepsilon'$ with
$0 < \varepsilon' < \varepsilon / (1+n) $ and choose $\delta >0$ for 
such an $\varepsilon'$.  
We choose sufficiently large $M_0$ such that for
every $M \ge M_0$, $f_M(\gamma_j(y),y)=f(\gamma_j(y),y)$ holds for every
$d(y,c) \ge \delta$ and for every $j$.  
We can conclude that for every $M \ge M_0$, 
\[
 |f_M(\gamma_j(y),y) - f(\gamma_j(y),y)| <  \varepsilon
\]
for all $y \in K$.  We have shown the following proposition.  

\begin{prop} If a right Hilbert C${}^*$-module $X$ is of $n$-branch 
 class and $n \ge 2$, $(u_k)_{k \in {\bf N}}$ as above is a basis 
for $X$.  \label{prop:branch_basis}
\end{prop}

We note that if $c$ is an accumulation point in $K$,
$(u_k)_{k \in \Lambda}$ is actually a countably infinite basis.  
\par
Since $\verb!#!C(\gamma)$ is finite,  we put $C(\gamma) 
=\{c_1,\dots,c_m \}$,  where $c_i \ne c_{i'}$ for $i\ne i'$.  
We take sufficiently small open neighborhoods $U_i$ of $c_i$ such that 
$C(\gamma) \cap \overline{U}_i = \{ c_i \}$ 
for $1 \le i \le m$ and $\overline{U}_i \cap \overline{U}_{i'} = \phi$ 
for $i \ne i'$.  We take an open neighborhood $V_i$ of $c_i$ such 
that $\overline{V}_i \subset U_i$ for each $i$.  We put $U_{m+1} = K
\backslash \bigcup_{i=1}^{m} \overline{V}_i$.  Then $\{U_i\}_{i=1}^{m+1}$
is an open covering of $K$.  
\par
We put $A_i={\rm C}(\overline{U}_i)$, 
$\C_i =\{(x,y) \in \C | y \in \overline{U}_i \} \subset \C$ and 
$X_i = {\rm C}(\C_i)$, $i=1,\cdots,m+1$.  
Then $X_i$'s are right Hilbert $A_i$ module naturally.  
We fix a $c_i \in C(\gamma)$.  
We put $\{ b^i_1,\dots,b^i_{r_i}\} = \gamma(c_i)$, 
where $b^i_s$ and $b^i_{s'}$ are different for $s \ne s'$.  
We put $\C_{i,s} = \{(x,y) \in   
\C_{i} | x=\gamma_k(y) \text{ for some } k \in I(b^i_s) \}$.  
Then we have $\C_{i} = \bigcup_{s=1}^{r_i}
\C_{i,s}$ and $\C_{i,s} \cap \C_{i,{s'}} = \phi$ for $s \ne s'$. 
We put $X_{i,s} = {\rm C}(\C_{i,s} )$.  
Then  $X_{i,s}$'s are right Hilbert $A_{i}$ module and we have 
\[
 (X_{i})_{A_{i}}
 = \bigoplus_{s=1}^{r_i} (X_{i,s})_{A_{i}}.  
\]
\par
Let $(u_{k}^{i,s})_{k \in \Lambda^{i,s}}$ be the basis for
$X_{i,s}$ defined in Proposition \ref{prop:branch_basis}.  
Put $\Lambda^i = \bigcup_{s=1}^{r_i}\{(s,k) | k \in 
\Lambda^{i,s})\}$.  

\begin{lemma}
 An indexed family $( u^{i,s}_k )_{ (s,k) \in \Lambda^i }$ 
is a basis for $X_{i}$ for each $i$.  
\label{lem:basis}
\end{lemma}
\begin{proof}
Since $X_i$ is a direct sum of  $X_{i,s}$'s, we can conclude the lemma by
 Lemma  \ref{lem:gather}.  
\end{proof}

Let $\{\psi_i\}_{i=1}^{m+1}$ be a partition of unity associated with
the open covering $\{U_i\}_{i=1}^{m+1}$.  
Let $( u^{i,s}_k )_{ (s,k) \in \Lambda^i }$ be a basis for
$X_i$ given by Lemma \ref{lem:basis}.  
We put 
$\tilde{u}^{i,s}_k (\gamma_j(y),y)=u^{i,s}_k(\gamma_j(y),y)
\psi_i(y)^{1/2}$.  
Then $\tilde{u}^{i,s}_k$'s can be extended to a function on $K$.  
Put $\Lambda = \bigcup_{i=1}^{m+1} \{ (i,(s,k)) | (s,k) 
\in \Lambda^i \}$.  

\begin{thm} Let $(K,d)$ be a compact metric space, 
$\gamma = (\gamma_1,\dots, \gamma_N)$ be a system of proper contractions
 on $K$.  We assume that $\gamma$ satisfies the finite branch 
 condition.  Then $( \tilde{u}^{i,s}_k)_{(i,(s,k)) \in \Lambda}$ 
defined above is a basis for $X$.  
 \label{thm:cover}
\end{thm}
\begin{proof}
 Let $f \in X$.  Since $\supp (f \cdot \psi_i) \subset U_i$, 
\[
 (f \cdot \psi_i)(\gamma_j(y),y) 
= \sum_{(k,s) \in \Lambda^i}
 u^{i,s}_k(u^{i,s}_k|f \cdot \psi_i)_A(\gamma_j(y),y).  
\]
uniformly.  
Using this equation we have 
\begin{align*}
  \sum_{i=1}^{m+1} \sum_{(s,k) \in \Lambda^i} \tilde{u}^{i,s}_k
    (\tilde{u}^{i,s}_{k}|f)_A (\gamma_j(y),y) 
 = & \sum_{i=1}^{m+1} \sum_{(s,k) \in \Lambda^i } u^{i,s}_{k} 
\cdot \psi_i^{1/2}
   (u^{i,s}_k \cdot \psi_i^{1/2}|f)_A (\gamma_j(y),y) \\
 = & \sum_{i=1}^{m+1} \sum_{(s,k) \in \Lambda^i} u^{i,s}_{k}
   (u^{i,s}_{k} |f \cdot \psi_i)(\gamma_j(y),y) \\
 = & \sum_{i=1}^{m+1} (f\cdot \psi_i)(\gamma_j(y),y) \\
 = & f(\gamma_j(y),y)
\end{align*}
uniformly with respect to $y$.  
\end{proof}

\begin{defn}
We call such basis as constructed in Theorem \ref{thm:cover} a patched 
 basis for $X$.   
\end{defn}

We put $\gamma_j^*(a)(y) =a(\gamma_j(y))$.  
For $a \in A$, we define a Borel function $\tilde{a}$ by 
\[
 \tilde{a}(y) = \sum_{x \in \gamma(y)}a(x)
 = \sum_{j=1}^N \frac{1}{e(\gamma_j(y),y)}a(\gamma_j(y)).  
\]
We note that if 
$C(\gamma)$ is not empty, $\tilde{a}$ is not continuous.

\begin{lemma}
 Let $X$ be an $n$-branch module.  Then for the basis 
$(u_k)_{k \in {\bf N}}$ constructed in 
Proposition \ref{prop:branch_basis} we have 
\[
 \sum_{k=1}^{\infty} (u_k|au_k)_A(y) = \tilde{a}(y)
\]
for every $y \in K$.  The left hand side converges unconditionally.
\label{lem:n_branch_sum}
\end{lemma}
\begin{proof}
 If $n=1$, we have 
\[
(u_1 | a u_1)_A(y)  = a(\gamma_1(y)).  
\]
We assume $n \ge 2$.  Then we have 

\begin{align*}
 \sum_{k=1}^{\infty}(u_k|au_k)_A(y) 
 = & \sum_{k=1}^{\infty} \sum_{j=1}^n \overline{u_k(\gamma_j(y),y)}
   a(\gamma_j(y))u_k(\gamma_j(y),y) \\
 = & \sum_{k=1}^{\infty} \sum_{j=1}^n |u_k(\gamma_j(y),y))|^2 
   a(\gamma_j(y))) \\
 = & \frac{1}{n} \sum_{j=1}^n a(\gamma_j(y)) +
  \frac{n-1}{n} \sum_{i=1}^{\infty} \sum_{j=1}^nv_i(d(\gamma_j(y),y))^2
 a(\gamma_j(y)). 
\end{align*}
The last expression is equal to $\sum_{j=1}^n a(\gamma_j(y))$ if 
$y \ne c$ and equal to $\frac{1}{n}\sum_{j=1}^n a(\gamma_j(y))$ if 
$y = c$.  
In any case, this is equal to $\tilde{a}(y)$.  If $a \in A^+$, the left
 hand side is monotone convergent, so we conclude that the left hand
 side converges unconditionally for general $a \in A$.  
\end{proof}

Since the left hand side in Lemma \ref{lem:n_branch_sum} converges
unconditionally, we may write 
\[
  \sum_{k \in \Lambda}(u_k|au_k)_A(y) = \tilde{a}(y).  
\]
\par
We take a base $(u^{i,s}_{k})_{k \in \Lambda^{i,s}}$ as 
in Theorem \ref{thm:cover}.  
Then 
$\sum_{s=1}^{r_i} \sum_{k \in \Lambda^{i,s}}
 (u^{i,s}_k|a u^{i,s}_k)_A(y)$ is $\sum_{s=1}^{r_i}a(b_s^i)$ for $y=c_i$ and 
$\sum_{s=1}^{r_i}(\sum_{j \in I(b_s^i)} a(\gamma_j(y))$ for $y \ne
s_i$.  This is equal to $\sum_{x \in \gamma(y)} a(x)$ in any case.  
We take a partition of unity $\{ \psi_i \}_{i=1}^{m+1}$ and
$(\tilde{u}^{i,s}_k)_{(i,(s,k)) \in \Lambda}$ as 
in Theorem \ref{thm:cover}.  We have 
\begin{align*}
 \sum_{i=1}^{m+1} \sum_{s=1}^{r_i} \sum_{k \in \Lambda^{i,s}}
 (\tilde{u}^{i,s}_k |a \tilde{u}^{i,s}_k)_A(y)
 & =  \sum_{i=1}^{m+1} \psi_i(y) \sum_{s=1}^{r_i}
   (u^{i,s}_k | a u^{i,s}_k)_A(y)(y) \\
 = &  \sum_{i=1}^{m+1} \psi_i(y)  \sum_{x \in \gamma(y)} a(x) \\
 = &  \sum_{x \in \gamma(y)} a(x) 
\end{align*}
Then we have the following proposition.  

\begin{prop} We assume that $\gamma$ satisfies the finite branch set
 condition.   For a patched basis $(u_k)_{k \in \Lambda}$ for
 $X_{\gamma}$, 
we have $\sum_{k \in \Lambda}(u_k|a u_k)_A(y) = \tilde{a}(y)$.  We note 
that the left side is monotone convergent for a positive $a \in A$.  
\label{prop:basis}
\end{prop}

Let $a \in I_X$.  Then $a(x)=0$ for $x \in B(\gamma)$.  Since $e(x,y)=1$ for $x
\ne B(\gamma)$, we have $\tilde{a}(y) = \sum_{j=1}^N a(\gamma_j(y))
= \sum_{j=1}^N \gamma_j^*(a)(y)$.  
\par
For a probability measure $\mu$ on $(K,d)$, $\tau^{\mu}$ denotes the
corresponding tracial state on $A$.  For a bounded Borel function $a$ on $K$, we
may define $\tau^{\mu} (a)$ by $\int_K a(y) d\mu (y)$.  

\begin{thm} Let $(K,d)$ be a self-similar set with respect to $\gamma$,  
 and $\gamma$ satisfy the finite branch condition.  
Let $\varphi$ be a $\beta$-KMS state on $\O_{\gamma}$.  Then the Borel
 probability measure $\mu$ on $(K,d)$ corresponding to the restriction 
of $\varphi $ to $A$ satisfies the  the following (3) and (4).  
\begin{align*}
  \sum_{j=1}^{N} \tau^{\mu}(\gamma_j^* (a)) & =  \lambda \tau^{\mu} (a)  
\quad (\forall a|_{B(\gamma)} = 0 ) \quad (3)  \\
  \tau^{\mu}(\tilde{a}) & \le  \lambda \tau^{\mu} (a) 
\quad (\forall a \in A^{+}) \quad (4), 
\end{align*}
where $\lambda = e^{\beta}$.  
\par
Let $\mu$ be a probability measure on $(K,d)$ satisfying (3) and (4).  
Then we can construct a $\beta$-KMS state $\varphi$ on $\O_{\gamma}$
 whose restriction to $A$ is $\tau^{\mu}$.  Moreover such a $\varphi$ is
 unique.  
\label{thm:KMS2}
\end{thm}
\begin{proof} 

Proposition \ref{prop:basis}, remark after Proposition \ref{prop:basis} 
and the monotone convergence theorem in integration theory show the 
theorem.  
\end{proof}

In the following, we present results which hold for general
situation.  

\begin{lemma}
 If a Borel probability measure $\mu$ on $K$ satisfies 
$\sum_{j=1}^N \tau^{\mu} (\gamma_j^*(a)) =  \lambda \tau^{\mu}(a)$ for 
arbitrary $a \in A$, a positive constant $\lambda$ must be equal to $N$.  
\label{lem:d1}
\end{lemma}
\begin{proof}
Putting $a=1$, we have $\lambda = N$.  
\end{proof}

\begin{lemma}(Hutchinson \cite{H})
Let $(K,d)$ be a compact metric space and $\gamma$ be a system of proper
 contraction.  Then there exists a unique measure $\mu$ on $(K,d)$ such that 
\[
\sum_{j=1}^N \tau^{\mu}(\gamma_j^*(a)) = N \tau^{\mu} (a)
\]
 for arbitrary $a \in A$. 
\label{lem:d2}
\end{lemma}

\begin{defn}
We denote by $\mu^H$ the measure given by Lemma \ref{lem:d2}, and call
 this measure the Hutchinson measure.  
\end{defn}

Then $\tau^{\mu^H}$ can be extended to a $\log N$-KMS sate 
$\varphi^H$ on $\O_{\gamma}$.  

\begin{rem} $\varphi^H$ is a KMS state of infinite
 type defined in Laca and Neshveyev \cite{LN}.  
When there exists no branched value 
for $\gamma$, $\tau^{\mu^H}$ is the unique KMS state on $\O_{\gamma}$.  
\end{rem}

\begin{defn}
 Let $\mu$ be a Borel probability measure on $(K,d)$.  We define 
 $c_{\mu}(x)$ by  $c_{\mu}(x) = \mu(\{x \})$.  
When $c_{\mu}(x)>0$,  we say that $\mu$ has a point mass at $x$.  
\end{defn}
We note that $c_{\mu}(x)$ is given by $\tau^{\mu}(\chi_{\{ x \}})$ 
where $\chi_{\{ x \}}$ is the characteristic function on a singleton
$\{x\}$.  

\begin{lemma}
 Let $x \in K$ and assume $\gamma^{-1}(x)=\{y_1,\dots,y_p\}$ with 
$y_i \ne y_j$  for $i \ne j$.  If $x \notin B(\gamma)$, we
 have $c_{\mu}(y_1) + \cdots + c_{\mu}(y_p) =\lambda c_{\mu}(x) $.  If 
$x \in B(\gamma)$, we have 
$c_{\mu}(y_1) + \cdots + c_{\mu}(y_p) \le \lambda c_{\mu}(x) $. 
\label{lem:equation}
\end{lemma}
\begin{proof}
Let $x \in K$.  We take $\{a_n\}_{n=1}^{\infty}$ such that $a_n \in A$
 and $a_n(x)$ tends to $\chi_{\{x\}}$ pointwise decreasingly.  
For each $j$,  we have $\gamma_j^*(a_n)$ tends to 
$\chi_{\{\gamma_j^{-1}(x)\}}$ pointwise decreasingly if $j \in I(x)$, 
and tend to zero pointwise decreasingly otherwise.  
We have $\tau^{\mu}(a_n) \to c_{\mu}(x)$, if $j \in I(x)$ 
$\tau^{\mu}(\gamma_j^{*}(a_n)) \to c_{\mu}(\gamma_j^{-1}(x))$ and 
otherwise $\tau^{\mu}(\gamma_j^{*}(a_n)) \to 0$.  
If $x \notin B(\gamma)$ we can take $a_n$ such that $a_n$ vanish on
 $B(\gamma)$ for all $n$.  
\begin{align*}
 \lambda c_{\mu}(x) & = \lambda \lim_{n \to \infty} 
 \tau^{\mu}(a_n)    = \lim_{n \to \infty} \sum_{j \in I(x)}
 \tau^{\mu}(\gamma_j^*(a_n))  \\
 &    = \sum_{j \in I(x)} c_{\mu}(\gamma_j^{-1}(x)) 
   = c_{\mu}(y_1) + \cdots + c_{\mu}(y_p).  
\end{align*}

If $x \in B(\gamma)$, we have 
\begin{align*}
 \lambda c_{\mu}(x)  & = \lambda \lim_{n\to \infty} \tau^{\mu}(a_n) 
   \ge \lim_{n \to \infty} \tau^{\mu}(\tilde{a_n})  \\
   & = \sum_{j \in I(x)} \frac{1}{e(x,\gamma_j^{-1}(x))}
            c_{\mu}(\gamma_j^{-1}(x))  \\
   & = c_{\mu}(y_1) + \cdots + c_{\mu}(y_p).    
\end{align*}

\end{proof}

\begin{lemma}
 Let $\mu$ be a Borel probability measure on $(K,d)$.  If  $\mu$ does
 not have a point mass at $B(\gamma)$ and satisfies (3) and (4) in
 Theorem \ref{thm:KMS2},  we can conclude that $\lambda = N$ and
 $\mu=\mu^H$.   
 If $\mu$ does not have a point mass at $B(\gamma) \cup 
\tilde{C}(\gamma)$ and satisfies (3), we get the same conclusion.  
\label{lem:total}
\end{lemma}
\begin{proof}
We assume that $\mu$ satisfies (4) and does not have a point mass at
 $B(\gamma)$.  We assume that $b \in B(\gamma)$ and
 $c_{1},\dots,c_{s}$ are mutually different elements in 
$\tilde{C}(\gamma)$ 
 such that $b = \gamma_{j_1}(c_{1})=\gamma_{j_2}(c_{2})= \cdots = 
\gamma_{j_s}(c_{s})$.  
Then by Lemma \ref{lem:equation} we have 
\[
  0 \le \sum_{p=1}^s c_{\mu}(c_{j_p}) \le \lambda c_{\mu}(b).  
\]
We conclude that $c_{\mu}(c_{j_p})=0$ for $1 \le p \le s$.  Since
 $c_{\mu}(x)=0$ for each $x \in B(\gamma)$, we have $c_{\mu}(y)=0$ for 
each $y \in \tilde{C}(\gamma)$.  
\par
We assume that $\mu$ satisfies (3) and does not have a point mass at
 $B(\gamma) \cup \tilde{C}(\gamma)$.  For each $a \in A^+$, there exists a 
 monotone increasing sequence of $\{ a_n \}_{n=1}^{\infty} \in A$ such 
that $a_n(x)=0$ for $x \in B(\gamma)$, and for $x \ne B(\gamma)$ 
there exits $n_0$ such that $a_n(x)=a(x)$ for $n \ge n_0$.  Then we have 
$a-a_n$ tends to $\sum_{x \in B(\gamma)} a(x) \chi_{\{ x \}}$ and 
$\gamma_j^*(a)-\gamma_j^*(a_n)$ tends to
 $\sum_{x \in (B(\gamma) \cap \gamma_j(K))}
 a(\gamma_j^{-1}(x)) \chi_{\{ \gamma_j^{-1}(x) \}}$.  
Since $\mu$ has no point mass on $B(\gamma) \cup \tilde{C}(\gamma)$, 
we have 
\[
 \lim_{n \to \infty} \tau^{\mu}(a_n)  = \tau^{\mu}(a)
 \qquad \lim_{n \to \infty} \tau^{\mu}(\gamma_j^*(a_n)) 
= \tau^{\mu}(\gamma_j^*(a)).  
\]
Then (3) holds for each $a \in A^+$.  By Lemma \ref{lem:d1} and Lemma
 \ref{lem:d2}, we have $\lambda = N$ and $\mu = \mu^H$.  
\end{proof}

\section{Classification of KMS states for specific examples} 
In this section, we present classifications of KMS states for some
specific examples.

\subsection{Dynamics on unit interval}
Let $K=[0,1]$, $d(x,y)=|x-y|$ and $\gamma=(\gamma_1,\dots,\gamma_N)$ 
be a system of proper contractions such that $K$ is self similar with
respect to $\gamma$.  
The following Lemma is easily verified.  
\begin{lemma}
We assume that $\gamma$ satisfies the open set condition.  
 For each $x \in B(\gamma)$ and $y \in C(\gamma)$ such that 
$x \in \gamma(y)$, 
we have $e(x,y)=2$.  Moreover $C(\gamma) =\tilde{C}(\gamma)$, and they are
 contained in $\{0,1\}$.  $B(\gamma)$ does not contain $0$ and $1$.  

\end{lemma} 
We always take $(0,1)$ as $V$.  We note that we have 
$B(\gamma) \subset V$.  We may assume that $\{\gamma_{j}\}_{j=1}^N$
satisfies 
$\gamma_1(1/2) < \gamma_2(1/2) < \cdots < \gamma_N(1/2)$.

In the following we assume that $\gamma$ satisfies the open set condition.
\begin{lemma}
 Let $y \in B(\gamma)$.  Then $O(y) \cap C(\gamma) = \phi$.  
 Let $y$ and $y'$ be distinct points in $B(\gamma)$.  Then
 $O(y) \cap O(y') = \phi$.  
\label{lem:orbit}
\end{lemma}

\begin{proof}
 Since $B(\gamma) \subset V$, $O(y)$ is contained in $V$.  We have
 $O(y) \cap \{ 0,1\} = \phi$.  
 We suppose that $\gamma_{i_1}\gamma_{i_2} \cdots \gamma_{i_n}(y) =
  \gamma_{j_1}\gamma_{j_2} \cdots \gamma_{j_m}(y')$ for 
$y$, $y' \in B(\gamma)$.  We assume that $n=m$.
 Then we have $i_1=j_1$, $i_2=j_2$ and $i_n=j_n$.  We have $y=y'$ and
 this is a contradiction.  We assume that
 $n \ge m+1$.  Then we have $y = \gamma_{n+1}\cdots \gamma_{m}(y')$, and 
this is a contradiction.  
\end{proof}

\begin{lemma}
 We assume that $\mu$ satisfies the conditions (3) and (4) in
 Theorem \ref{thm:KMS2}.
 When $\lambda > 1$, we have $c_{\mu}(0)=0$ and $c_{\mu}(1)=0$.  
 \label{lem:zero1}
\end{lemma}

\begin{proof}
We suppose that $\gamma_1(0)=0$ and $\gamma_N(0)=1$.  
Since we have $\gamma^{-1}(0)=\gamma_1^{-1}(0)=\{ 0 \}$, we have
 $c_{\mu}(0)=\lambda c_{\mu}(0)$ by Lemma \ref{lem:equation}.  
 Since $\lambda >1$ we have $c_{\mu}(0)=0$.  
Since we have $\gamma^{-1}(1)=\gamma_N^{-1}(1) = \{ 0 \}$, we have 
$c_{\mu}(0) = \lambda c_{\mu}(1)$ by Lemma \ref{lem:equation}.  
Then $c_{\mu}(1)=0$ follows.  
 \par
We suppose $\gamma_1(1)=0$ and $\gamma_N(0)=1$.  Then by a similar
 computation, we have $c_{\mu}(1)=\lambda c_{\mu}(0)$ and 
$c_{\mu}(0) = \lambda c_{\mu}(1)$.  Since $\lambda > 1$, 
we have $c_{\mu}(0)=0$ and $c_{\mu}(1)=0$.  
\par
For other two cases, we can prove lemma similarly.  
\end{proof}

\begin{lemma}
 We assume that $\mu$ satisfies (3) and (4) in Theorem \ref{thm:KMS2}
 and has a point mass at some point $y \in B(\gamma)$.  
 Then we have $\lambda > N$ and 
 $c_{\mu}(\gamma_{j_1}\cdots \gamma_{j_n}(y)) = \lambda^{-n} c_{\mu}(y)$
 for each $(j_1,\dots,j_n) \in \{1,\dots,N \}^n$.  
 \label{lem:coeff}
\end{lemma}

\begin{proof}
If $x=\gamma_j(y)=\gamma_{j'}(y')$ for $y$, $y' \in V$ and $1 \le j$, 
$j'\le N$, we have $y=y'$ and $j=j'$, and for $y \in B(\gamma)$
 $\gamma_{j_1}\cdots \gamma_{j_n}(y)$ is contained in 
$V \cap B(\gamma)^c$ for every $(j_1,\dots,j_n)$, $n \ge 1$.  
Then by Lemma \ref{lem:equation},  we have 
\[
 c_{\mu}(\gamma_{j_n}\cdots \gamma_{j_1}(y)) 
 = \lambda c_{\mu} (\gamma_{j_{n+1}}\gamma_{j_n}\cdots \gamma_{j_1}(y)).  
\]
for $n \ge 0$.  Then we have 
\[
 c_{\mu}(\gamma_{j_n}\cdots \gamma_{j_1}(y)) 
= \lambda^{-n} 
 c_{\mu}(y)
\]
By Lemma \ref{lem:orbit} $O(y) \cap C(\gamma)= \phi$.  Then 
\[
 \mu(1) \ge \sum_{n=1}^{\infty}\sum_{(j_1,\dots,j_n) 
\in \{1,\cdots,N\}^n}c_{\mu} 
 (\gamma_{j_1}\cdots \gamma_{j_n}(y))  = \sum_{n=1}^{\infty} 
 \left(\frac{N}{\lambda}\right)^n c_{\mu}(y).  
\]
This shows that $N < \lambda$ is necessary for $\mu$ to be bounded.  
\end{proof}

\begin{lemma}
 We assume that $\mu$ is a Borel probability measure on $[0,1]$ and
 satisfies (3) and (4) in Theorem \ref{thm:KMS2}.
 Then we have $\lambda \ge N$ and $c_{\mu}(0)=c_{\mu}(1) = 0$.  
 \label{lem:zero2}
\end{lemma}
\begin{proof} 
If $\mu$ does not have a point mass at $B(\gamma)$, by 
 Lemma \ref{lem:total} we have $c_{\mu}(y)=0$ for
 $y \in C(\gamma)=\tilde{C}(\gamma)$, and have $\mu = \mu^H$ 
 and $\lambda = N$.  If $\mu$ has a point mass at $B(\gamma)$,
 then by Lemma \ref{lem:coeff} we have
 $\lambda > N$.  By Lemma \ref{lem:zero2},  
we have $c_{\mu}(0)=c_{\mu}(1)=0$.  
 \end{proof}

\begin{lemma} 
 A probability measure $\mu$ on $[0,1]$ satisfying (3) and (4) in
 Theorem \ref{thm:KMS2} can have 
point mass only at $\{ O(y) | y \in B(\gamma)\}$.
 In particular, $\mu^H$ has no point mass.  
\end{lemma}
\begin{proof} 
Let $ y \notin \{ O(y) | y \in B(\gamma) \}$ and $y \ne 0$, $1$.  
Then we can construct a sequence $\{y_i \}_{i=0}^{\infty}$ such that 
$y_0=y$ and $y_i = \gamma_{j_i}(y_{i+1})$.  There exist three 
possibilities.  
\begin{enumerate}
 \item All $\{y_0,y_1,\cdots,y_i,\cdots \}$ are different 
and $y_i \in V \cap B(\gamma)^c$ for all $i$.  
  \item There exists $i_0$ and $m \ge 1$ such that $y_0$, $\cdots$
       $y_{i_0}$ are all different, and $y_{i_0 + m}=y_{i_0}$ and 
   $y_i \in V \cap B(\gamma)^c$ for all $i$.  
 \item There exists $i_0$ such that $y_0$, $\cdots$ $y_{i_0-1}$ are all
       different and in $V \cap B(\gamma)^c$ and $y_{i_0}=0$ or $1$.  
\end{enumerate}
In case (1), we have $c_{\mu}(y_{i+1}) = \lambda c_{\mu}(y_i)$.  Then we
 have $c_{\mu}(y_i) = \lambda^i c_{\mu}(y)$.  This shows that if 
 $c_{\mu}(y) > 0$ $c_{\mu}(y_i) \to \infty$, and is a contradiction.  
 In case (2), we have $c_{\mu}(y_{i_0}) = c_{\mu}(y_{i_0 + m}) =
 \lambda^{m}c_{\mu}(y_{i_0})$.  Since $\lambda > 1$, we have
 $c_{\mu}(y_0) = 0$.  $c_{\mu}(y) = \lambda^{- i_0 } c_{\mu}(y_{i_0})$ 
shows that $c_{\mu}(y)=0$.  
 In case (3), we have $c_{\mu}(0)+  c_{\mu}(1) =
 \lambda c_{\mu}(y_{i_0-1})$, $c_{\mu}(0) = \lambda c_{\mu}(y_{i_0-1})$ or 
$c_{\mu}(1) = \lambda c_{\mu}(y_{i_0-1})$.  This shows that
 $c_{\mu}(y_{i_0-1})=0$.  Then we have $c_{\mu}(y) =
 \lambda^{-i_0} c_{\mu}(y_{i_0 -1}) =0$.  
In any case, we can conclude $c_{\mu}(y)=0$.  
\end{proof}

Let $y \in B(\gamma)$.  We consider 
 \[
 a \in A^+ \to \sum_{n=0}^{\infty} \sum_{(j_1,\cdots,j_n)\in \{1,\dots,N\}^n}
 ( 1 / \lambda )^n a (\gamma_{j_1}\cdots \gamma_{j_n}(y)).  
 \]
This gives a bounded Borel measure on $[0,1]$ if and only if $\lambda >
N$.  When $\lambda > N$, we can define a Borel probability measure 
$\mu_{y,\lambda}$ by 
\[
  \mu_{y,\lambda} = \frac{\lambda - N}{\lambda}
  \sum_{n=0}^{\infty} \sum_{(j_1,\cdots,j_n) \in \{1,\dots,N\}^n}
 ( 1 / \lambda )^n \delta_{\gamma_{j_1}\cdots \gamma_{j_n}(y)}
 \]
for $a \in A$.  
We put $\beta = \log \lambda$.  
\begin{prop}
 The measure $\mu_{y,\lambda}$ satisfies the condition (3) and (4) in
 Theorem \ref{thm:KMS2}, and is extended to a $\beta$-KMS state 
$\varphi^{y,\lambda}$ on $\O_{\gamma}$.  
\end{prop}
\begin{proof}
 Let $y \in B(\gamma)$.  Since 
\[
 \tau^{\mu_{y,\lambda}}(a) = \frac{\lambda -N}{\lambda}
 \sum_{n=0}^{\infty}\sum_{(j_1,\cdots,j_n)
 \in {\{1,\cdots,N \}}^n } \lambda^{-n} a(\gamma_{j_1}\cdots 
\gamma_{j_n}(y)), 
\]
we have 
\[
 \lambda \tau^{\mu_{y,\lambda}}(a) 
=  \frac{\lambda - N}{\lambda} \sum_{n=0}^{\infty}
   \sum_{(j_1,\cdots,j_n) \in \{1,\dots,N \}^n} \lambda^{-n+1}
   a(\gamma_{j_1}\cdots  \gamma_{j_n}(y)).  
\]
Since by Lemma \ref{lem:orbit} $\gamma_{j_1} \cdots \gamma_{j_n}(y)$ 
are not contained in $C(\gamma)$ for every $n \ge 0$,  we have 
\begin{align*}
     \tau^{\mu_{y,\lambda}}(\tilde{a}) 
  = &  \frac{\lambda - N}{\lambda} \sum_{n=0}^{\infty} 
  \sum_{(j_1,\cdots,j_n) \in {\{1,\cdots,N\}}^n} \lambda^{-n} 
  \tilde{a}(\gamma_{j_1}\cdots \gamma_{j_n}(y)) \\
  = &  \frac{\lambda - N}{\lambda} \sum_{j=1}^N  \sum_{n=0}^{\infty}
   \sum_{(j_1,\cdots,j_n) \in {\{1,\cdots,N\}}^n} \lambda^{-n} 
   a(\gamma_j \gamma_{j_1}\cdots \gamma_{j_n}(y)) \\
   = &  \frac{\lambda - N}{\lambda} \sum_{n=1}^{\infty}
  \sum_{(\tilde{j}_1,\cdots,\tilde{j}_n) \in {\{1,\cdots,N \}}^n} \lambda^{-(n-1)} 
   a(\gamma_{\tilde{j}_1}\cdots \gamma_{\tilde{j}_n}(y)).  
\end{align*} 
We have 
\[
 \lambda \tau^{\mu_{y,\lambda}}(a) = \tau^{\mu_{y,\lambda} }(\tilde{a}) 
   + \left( \frac{\lambda -N}{\lambda} \right) a(y).  
\]
If $a \in I_X$, we have $\lambda \tau^{\mu_{y,\lambda}}(a) 
=\tau^{\mu_{y,\lambda}}(\tilde{a})$ because $a(y)=0$.  
If $a$ vanish on $B(\gamma)$, we have $a(y) \ge 0$ and have
 $\tau^{\mu_{y,\lambda}}(\tilde{a}) \le \lambda 
\tau^{\mu_{y,\lambda}}(a) $.  
\end{proof}


\begin{rem} $\varphi^{y,\lambda}$ is a KMS state of finite type defied
 in Laca and Neshveyev \cite{LN}.  
\end{rem}

\begin{thm}
 Let $\gamma$ be a system of proper contractions on $[0,1]$ and satisfy
 the open set condition.  Then a $\beta$-KMS state on $\O_{\gamma}$ with
 respect to the gauge action exists only if $\lambda = e^{\beta} \ge N $
 and are classified as follows:
\begin{enumerate}
 \item When $\lambda = N$, $\varphi^H$ is the unique KMS state.  
 \item When $\lambda > N$, $\beta$-KMS state is expressed by
       a convex combination of $\{ \varphi^{y,\lambda} |
       y \in B(\gamma) \}$.  
\end{enumerate}
Moreover $\varphi^H$ is unique $\log N$-KMS state, and if
 $B(\gamma)$ is not empty, $\varphi_{y,\lambda}$ is an extreme 
$\log \lambda$-KMS state.  
\end{thm}
\begin{proof} 
 Let $\mu$ be a Borel probability measure on $[0,1]$ and satisfy the
 condition (3) and (4) in Theorem \ref{thm:KMS2}.
 Then by Lemma \ref{lem:zero2}, we have $\mu \ge N$.  
 If $\lambda = N$ then $\mu$ dose not have point mass at $B(\gamma)$,
 and then by Lemma \ref{lem:total} we have $\mu = \mu^H$.  
\par
We assume that $\lambda > N $.  By Lemma  \ref{lem:coeff} and 
Lemma \ref{lem:zero1}, $\mu - \sum_{y \in B(\gamma)} 
 c_{\mu}(y) \mu_{y,\lambda}$ is a positive Borel measure, satisfies the
 condition (3) and does not have point mass at $B(\gamma) \cup
 \tilde{C}(\gamma)$.  Then by Lemma \ref{lem:total}, the condition (3)
 must hold for  all $a \in A$.  
 Then we have $\mu - \sum_{y \in B(\gamma)} c_{\mu}(y) \mu_{y,\lambda}=0$.  
\par
Lastly, we show that $\varphi^{y,\lambda}$ is extreme.  
We write $\varphi^{y,\lambda} = t \varphi_1 + (1-t) \varphi_2 $, where
 $\varphi_1$ and $\varphi_2$ be a $\beta$-KMS state on $\O_{\gamma}$ and
 $0 < t < 1$.  By restricting to $A$, we conclude 
that $\varphi_i = \varphi^{y,\lambda } $.  This shows that
 $\varphi^{y,\lambda}$ is extreme.  
\end{proof}

We assume that $\gamma$ is a section of a expansive map $\gamma^{-1}$
on $K$.  
The value of inverse temperature of KMS states have a relation with 
entropy of $\gamma^{-1}$.  

\begin{prop}
 The minimum value of the logarithm of the inverse temperature 
 of KMS states on $\O_{\gamma}$ is equal to the entropy of the 
 map $\gamma^{-1}$ on $K$.  
\end{prop}

\begin{proof}
 By Theorem 7.2 in \cite{MS}, the entropy $h(\gamma^{-1})$ is equal 
 to $\log N$.  This is equal to the minimum value of the logarithm of 
the inverse temperature of KMS states on $\O_{\gamma}$.  
\end{proof}

 We consider the example \ref{exam:tent} i.e. 
 $K=[0,1]$, $\gamma_1(y) = \frac{1}{2}y$,  
and $\gamma_2(y) = 1 - \frac{1}{2}y$.  We denote by $\O_{\text{tent}}$ 
the Cuntz-Pimsner C${}^*$-algebra for this example.  
We note that 
$B(\gamma) =\{ \frac{1}{2}\}$ and  $C(\gamma) = \{1\}$.  
Let $\mu$ be the normalized Lebesgue measure on $[0,1]$.  
We assume $\lambda >2$.  We put 
\[
\mu_{1/2,\lambda} = \frac{\lambda -2}{\lambda} 
        \sum_{n=0}^{\infty}\sum_{(j_1,\cdots,j_n) \in \{ 1,2 \}^n} 
	 \lambda^{-n} \delta_{\gamma_{j_1}\cdots \gamma_{j_n}(1/2)} 
\]
Then we have the following: 

\begin{prop}
A $\beta$-KMS state on $\O_{\text {tent}}$ exist if and only if 
$\beta=\log \lambda \ge \log 2$.  If $\lambda = 2$, $\beta$-KMS 
state is unique and given by $\varphi^{\mu}$.  
If $\lambda >2$, $\beta$-KMS state is unique and given by 
$\varphi^{1/2,\lambda}$.   
\end{prop}

 We consider example \ref{exam:Cuntz} ie $K = [0,1]$, 
$\gamma_1(y)=\frac{1}{2}y$ and $\gamma_2(y) = \frac{1}{2}(y+1)$. 
In this case, $B(\gamma)=\phi$ 
and $C(\gamma)=\phi$.  $\beta$-KMS state on $\O_{\gamma}$ exists if and 
only if $\lambda = \log 2$ and given by the normalized Lebesgue measure 
on $[0,1]$.

\subsection{Sierpinski Gasket}
As the case of dynamics on unit interval, we can classify KMS states for 
the C${}^*$-algebra associated with Sierpinski Gasket introduced in 
Kajiwara-Watatani \cite{KW2}.  The contractions in this example 
are considered to be cross sections for rational map on Riemaniann 
sphere whose Julia set is homeomorphic to Sierpinski Gasket.  
\par
Let $\Omega$ be a regular triangle in $\bR^2$ with three vertexes  
$c_1=(1/2,\sqrt{\,3}/2)$, $c_2 = (0,0)$ and $c_3=(1,0)$.  
The middle point of $c_1c_2$ is denote by $b_1$, the middle point of 
$c_1c_3$ is denoted by $b_2$ and the middle point of $c_2c_3$ is 
denoted by $b_3$.  
We define proper contractions $\tilde{\gamma_i}$ $(i=1,2,3)$ by 
\[
 \tilde{\gamma}_1(x,y)  = \left( \frac{x}{2} + \frac{1}{4}, \frac{y}{2} 
+ \frac{\sqrt{\,3}}{4}\right), \quad
 \tilde{\gamma}_2(x,y)  = \left( \frac{x}{2}, \frac{y}{2}\right), \quad 
 \tilde{\gamma}_3(x,y)   = \left( \frac{x}{2} + \frac{1}{2},\frac{y}{2} 
 \right).  
\]
Let $\alpha_{\theta}$ be a rotation by the angle $\theta$.  
We put $\gamma_1 = \tilde{\gamma}_1$, $\gamma_2 = \alpha_{-2\pi/3} \circ  
\tilde{\gamma}_2$ and $\gamma_3 = \alpha_{2 \pi/3} \circ 
\tilde{\gamma}_3$.   We denote by $S$ with the metric $d$ induce from 
${\bf R}^2$ the self similar set determined by
$\gamma=(\gamma_1,\gamma_2,\gamma_3 )$.    
We note that $c_i$ and $b_i$ $i=1,2,3$ are contained in $S$.  
Putting $V = S \backslash \{c_1,c_2,c_3\}$, $\gamma$ satisfies the open 
set condition.  
In this case, we have $B(\gamma)=\{b_1,b_2,b_3\}$ and $C(\gamma) 
=\tilde{C}(\gamma) =\{ c_1,c_2,c_3 \}$, and $\gamma$ satisfies the 
finite branch condition.  
We denote by $\O_{S,\gamma}$ Cuntz-Pimsner algebra constructed 
from $S$ and the above $\gamma$.  
\par
Let $\mu$ be the Borel probability measure on $(S,d)$ satisfying the 
condition (3) and (4).  
We get the conditions of point mass of $\mu$.  
\begin{lemma}
If $\lambda > 1$, we have $c_{\mu}(c_1)=0$, $c_{\mu}(c_2)=0$ and 
 $c_{\mu}(c_3)=0$.  
\end{lemma}

\begin{proof}
 We note that $\gamma^{-1}(c_1) = \{c_1\}$, $\gamma^{-1}(c_2) = \{c_3\}$ 
and $\gamma^{-1}(c_3) = \{c_2\}$.  By Lemma \ref{lem:equation}, we have 
\[
 c_{\mu}(c_1) = \lambda c_{\mu}(c_1) \qquad 
 c_{\mu}(c_2) = \lambda c_{\mu}(c_3) \qquad 
 c_{\mu}(c_3) = \lambda c_{\mu}(c_2).  
\]
When $\lambda > 1$, then these show that $c_{\mu}(c_1)=0$, 
 $c_{\mu}(c_2)=0$ and $c_{\mu}(c_3)=0$.  
\end{proof}

As in the case of dynamics on unit interval, we have the following 
Lemmas.  

\begin{lemma} 
For $y \in B(\gamma)$, we have $O(y) \cap C(\gamma)= \phi$, and for $y$, 
 $y'\in B(\gamma)$ with $y \ne y'$, we have $O(y) \cap O(y') = \phi$.  
\end{lemma}

\begin{lemma}
If $\mu$ satisfies the condition (3) and (4) in Theorem \ref{thm:KMS2} 
 and does not have a point mass at $B(\gamma)$, or if $\mu$ satisfies
 the condition (3) in Theorem \ref{thm:KMS2} and does not have a point
 mass at $B(\gamma) \cup C(\gamma)$, we have $\lambda=3$ 
 and $\mu = \mu_H$.  
\end{lemma}

\begin{lemma}
If $\mu$ satisfying the condition (3) and (4) in
 Theorem \ref{thm:KMS2} has a point mass at $B(\gamma)$, then
 we have $\lambda >3$ and $c_{\mu}(\gamma_{j_1}\cdots \gamma_{j_n}(y)) = \lambda^{-n}c_{\mu}(y)$.  
\end{lemma}

\begin{lemma} If $\mu$ has the condition (3) and (4) in
 Theorem \ref{thm:KMS2}, then $\mu$ does not have a point mass at 
$C(\gamma)$. 
\end{lemma}

Let $\lambda > 3$.  As in dynamics on unit interval, for $y \in 
B(\gamma)$, then we define a probability measure $\mu_{\lambda,y}$ as
follows: 
\[
 \tau^{\mu_{y,\lambda}}(a) = \frac{\lambda-3}{\lambda} 
 \sum_{n=0}^{\infty}\sum_{(j_1,\dots,j_n) 
 \in \{1,2,3\}^n} a(\gamma_{j_1} \cdots \gamma_{j_n}(y)).  
\]

As dynamic for unit interval, we have the following: 

\begin{lemma}
$\mu_{\lambda,y}$ satisfies the condition (3) and (4), and is extended 
 to the $\log \lambda$-KMS state on $\O_{S,\gamma}$.  
\end{lemma}

We can get classification of KMS states on $\O_{S,\gamma}$.  
Let $\beta = \log \lambda$.  

\begin{thm}  Let $S$ be the Sierpinski gasket defined by contractions 
$\gamma$ as above. Then $\beta$-KMS state on $\O_{S,\gamma}$ with 
 respect to the gauge action exists only if $\lambda \ge 3$ 
 and are classified as follows:  
\begin{enumerate}
 \item When $\lambda = 3$, $\varphi^H$ is the unique KMS state.  
 \item When $\lambda > 3$, each $\beta$-KMS state is expressed by 
       a convex combination of  \par 
       $\{ \varphi^{y,\lambda}\, | \,y=b_1,\, b_2,\,b_3 \}$ 
\end{enumerate}
Moreover $\varphi_{y,\lambda}$'s are an extreme $\log \lambda$-KMS state.   
\end{thm}


%


\begin{thebibliography}{99}
\bibitem{BR} Bratteli O. and Robinson D.W., 
{\it Operator algebras and quantum statistical mechanics II.  
Equilibrium states.  Models in quantum statistical mechanics.,}
Texts and Monographs in Physics. Springer-Verlag, Berlin, 1997


\bibitem{EFW} Enomoto M., Fujii M. and Watatani Y., 
{\it KMS states for gauge action on $O_A$, }
Math. Japon. 29(1984), 607--619


\bibitem{Ev} Evans D.,
{\it On O${}_n$, } Publ. Res. Inst. Math. Sci. Kyoto Univ.
16(1980)


\bibitem{Ex1} Exel R.,
{\it Crossed-products by finite index endomorphisms and KMS states,}
J. Funct. Anal. 199(2003), 153--183


\bibitem{Ex2} Exel R., 
{\it KMS states for generalized gauge actions on Cuntz-Krieger
algebras (An application of the Ruelle-Perron-Frobenius Theorem)}
[arXiv : math.OA/0110183]


\bibitem{EL} Exel R. and Laca M.,  
{\it Partial dynamical systems and the KMS condition,} 
Commun. Math. Phys., 232(2003), 223--277


\bibitem{FMR} Fowler N.J., Muhly P.S. and Raeburn I., 
{\it Representations of Cuntz-Pimsner Algebras},  
Indiana Univ. Math. J. {\bf 52} (2003), 569--605.


\bibitem{H} Huchchinson J.E., 
{\it Fractals and self-similarity,} Indiana Univ. Math. J. 30 (1981), 713--747.


\bibitem{KPW} Kajiwara T., Pinzari C. and Watatani Y.,
{\it Jones index theory for Hilbert C${}^*$--bimodules and its equivalence
with conjugation theory, }
J. Funct. Anal. in press


\bibitem {KW1} Kajiwara T. and Watatani Y., 
{\it C${}^*$-algebra associated with complex dynamical systems, }
[arXiv : math.OA/0309293 ]


\bibitem {KW2} Kajiwara T. and Watatani Y., 
{\it C${}^*$-algebras associated with self-similar sets}
[arXiv : math.OA/0312481 ] 


\bibitem{KP} Kerr D. and Pinzari C.,
{\it  Noncommutative pressure and the variational principle in 
Cuntz-Krieger-type C${}^*$-algebras. }
J. Funct. Anal. 188(2002), 156--215


\bibitem{KR} Kumjian A. and Renalut J., 
{\it KMS states on C${}^*$-algebras associated to expansive maps, }
[arXiv : math. OA/0305044]


\bibitem{LN} Laca M. and Neshveyev S.,
{\it KMS states of quasi-free dynamics on Pimsner algebras,}
J. Funct. Anal. 211(2004) 457--482


\bibitem{MWY} Matsumoto K., Watatani Y. and Yoshida M.
{\it KMS-states for gauge action on C${}^*$-algebras associated with
	subshifts,}
Math. Z. 228(1998), 489--509


\bibitem{MS} Melo W. and Strien S., 
{\it One-dimesional Dynmics, }
Springer,  1993


\bibitem{MT} Muhly P.S. and Tomforde M.,
{\it Topological quivers,}
[arXiv :  math.OA/0312109]



\bibitem{OP} Olsen D. and Pederson G.K., 
{\it Some C${}^*$-dynamical systems with a single KMS state,}
Math. Scand. 42(1978), 111--118


\bibitem{Pi} Pimsner M., 
{\it A class of $C^*$-algebras generating both Cuntz-Krieger algebras
and crossed product by ${\mathbb Z}$,} 
Free probability theory, AMS, (1997), 189--212.


\bibitem{PWY} Pinzari C., Watatani Y. and Yonetani K.
{\it KMS states, entropy and the variational principle in full 
C${^*}$ -dynamical systems,} 
Commun. Math. Phys. 213(2000), 331--379


\bibitem{R1} Renault, J., 
{\it A groupoid approach to C${}^*$-algebras}
Lecture Notes in Mathematics, 793. Springer, Berlin, 1980.


\bibitem{R2} Renault, J., 
{\it Cuntz-like algebras, }
Operator theoretical methods (Timicsoara, 1998), (2000), 371--386,
Theta Found., Bucharest


\end{thebibliography}
\end{document}